\documentclass[12pt]{amsart}
\usepackage{amssymb,mathrsfs}

\begin{document}


\title{Invariant differential operators
on the Minkowski-Euclid space  }

\author{Jae-Hyun Yang}
\address{Department of Mathematics, Inha University, Incheon 402-751, Korea}
\email{jhyang@inha.ac.kr }


\newtheorem{theorem}{Theorem}[section]
\newtheorem{lemma}{Lemma}[section]
\newtheorem{proposition}{Proposition}[section]
\newtheorem{remark}{Remark}[section]
\newtheorem{definition}{Definition}[section]

\renewcommand{\theequation}{\thesection.\arabic{equation}}
\renewcommand{\thetheorem}{\thesection.\arabic{theorem}}
\renewcommand{\thelemma}{\thesection.\arabic{lemma}}
\newcommand{\bbr}{\mathbb R}
\newcommand{\bbs}{\mathbb S}
\newcommand{\bn}{\bf n}
\newcommand\charf {\mbox{{\text 1}\kern-.24em {\text l}}}
\newcommand\fg{{\mathfrak g}}
\newcommand\fk{{\mathfrak k}}
\newcommand\fp{{\mathfrak p}}
\newcommand\g{\gamma}
\newcommand\G{\Gamma}
\newcommand\ka{\kappa}
\newcommand\al{\alpha}
\newcommand\be{\beta}
\newcommand\lrt{\longrightarrow}
\newcommand\s{\sigma}
\newcommand\ba{\backslash}
\newcommand\lmt{\longmapsto}
\newcommand\CP{{\mathcal P}_n}
\newcommand\CM{{\mathcal M}}
\newcommand\BC{\mathbb C}
\newcommand\BZ{\mathbb Z}
\newcommand\BR{\Bbb R}
\newcommand\BQ{\mathbb Q}
\newcommand\Rmn{{\mathbb R}^{(m,n)}}
\newcommand\PR{{\mathcal P}_n\times {\mathbb R}^{(m,n)}}
\newcommand\Gnm{GL_{n,m}}
\newcommand\Gnz{GL_{n,m}({\mathbb Z})}
\newcommand\Gjnm{Sp_{n,m}}
\newcommand\Gnml{GL(n,{\mathbb R})\ltimes {\mathbb R}^{(m,n)}}
\newcommand\Snm{SL_{n,m}}
\newcommand\Snz{SL_{n,m}({\mathbb Z})}
\newcommand\Snml{SL(n,{\mathbb R})\ltimes {\mathbb R}^{(m,n)}}
\newcommand\Snzl{SL(n,{\mathbb Z})\ltimes {\mathbb Z}^{(m,n)}}
\newcommand\la{\lambda}
\newcommand\GZ{GL(n,{\mathbb Z})\ltimes {\mathbb Z}^{(m,n)}}
\newcommand\DPR{{\mathbb D}(\PR)}
\newcommand\Rnn{{\mathbb R}^{(n,n)}}
\newcommand\Yd{{{\partial}\over {\partial Y}}}
\newcommand\Vd{{{\partial}\over {\partial V}}}
\newcommand\Ys{Y^{\ast}}
\newcommand\Vs{V^{\ast}}
\newcommand\DGR{{\mathbb D}(\Gnm)}
\newcommand\DKR{{\mathbb D}_K(\Gnm)}
\newcommand\DKS{{\mathbb D}_{K_0}(\Snm)}
\newcommand\fa{{\frak a}}
\newcommand\fac{{\frak a}_c^{\ast}}
\newcommand\SPR{S{\mathcal P}_n\times \Rmn}
\newcommand\DSPR{{\mathbb D}(\SPR)}
\newcommand\BD{{\mathbb D}}

\thanks{\endgraf Keywords and phrases\,: invariants, invariant differential
operators, the Minkowski-Euclid space.\\
\indent This work was supported by Inha University Research Grant.  }

\begin{abstract} For two positive integers $m$ and $n$, we let ${\mathcal
P}_n$ be the open convex cone in ${\mathbb R}^{n(n+1)/2}$
consisting of positive definite $n\times n$ real symmetric
matrices and let $\BR^{(m,n)}$ be the set of all $m\times n$ real
matrices. In this article, we investigate differential operators
on the non-reductive homogeneous space ${\mathcal P}_n\times {\mathbb
R}^{(m,n)}$ that are invariant under the natural action of the
semidirect product group $GL(n,\BR)\ltimes \BR^{(m,n)}$ on the Minkowski-Euclid space
${\mathcal P}_n\times {\mathbb R}^{(m,n)}$. These invariant differential operators
play an important role in the theory of automorphic forms on $GL(n,\BR)\ltimes \BR^{(m,n)}$
generalizing that of automorphic forms on $GL(n,\BR).$
\end{abstract}
\maketitle



\begin{section}{{\bf Introduction}}
\setcounter{equation}{0}

We let
$$\CP=\left\{\, Y\in \BR^{(n,n)}\,\vert\ Y=\,^tY>0\ \right\}$$ be
the open convex cone of positive definite symmetric real matrices of degree $n$
in the Euclidean space $\BR^{n(n+1)/2},$ where $F^{(k,l)}$ denotes
the set of all $k\times l$ matrices with entries in a commutative
ring $F$ for two positive integers $k$ and $l$ and $^t\!M$ denotes
the transpose matrix of a matrix $M$. Then the general linear
group $GL(n,\BR)$ acts on $\CP$ transtively by
\begin{equation}
g\cdot Y=gY\,^tg,\ \ \ \ g\in GL(n,\BR),\ Y\in \CP.
\end{equation}
Therefore $\CP$ is a symmetric space which is diffeomorphic to the
quotient space $GL(n,\BR)/O(n)$, where $O(n)$ denotes the
orthogonal group of degree $n$. A. Selberg [9] investigated
differential operators on $\CP$ invariant under the action (1.1)
of $GL(n,\BR)$ (cf. [6,\,7]).

\vskip 0.2cm We let
\begin{equation*}
\Gnm=GL(n,\BR)\ltimes \Rmn \end{equation*} be the semidirect
product of $GL(n,\BR)$ and the abelian additive group $\Rmn$
equipped with the following multiplication law
\begin{equation*}
(g,\lambda)\cdot (h,\mu)=\big(gh,\lambda\,^th^{-1}+\mu\big),
\end{equation*} where $g,h\in GL(n,\BR)$ and $\lambda,\mu\in \Rmn$.
Then we have the {\it natural action} of $\Gnm$ on the non-reductive
homogeneous space $\PR$ given by
\begin{equation}
(g,\lambda)\cdot (Y,V)=\big(\,
gY\,^tg,\,(V+\lambda)\,^tg\big),\end{equation} where $g\in
GL(n,\BR),\ \lambda\in \Rmn,\ Y\in \CP$ and $V\in \Rmn.$

\vskip 0.2cm For brevity, we set ${\mathcal P}_{n,m}=\PR$ and
$K=O(n).$ Since the action (1.2) of $GL_{n,m}$ is transitive,
${\mathcal P}_{n,m}$ is diffeomorphic to $\Gnm/K$. We observe that
the action (1.2) of $\Gnm$ generalizes the action (1.1) of
$GL(n,\BR).$

\vskip 0.2cm The reason why we study the non-reductive homogeneous space
${\mathcal P}_{n,m}$ may be explained as follows. Let
$$\Gamma_{n,m}=GL(n,\BZ)\ltimes \BZ^{(m,n)}$$
\noindent be the arithmetic subgroup of $GL_{n,m}$, where $\BZ$ is
the ring of integers. The arithmetic quotient
$\Gamma_{n,m}\backslash {\mathcal P}_{n,m}$ may be regarded as the
universal family of principally polarized real tori of dimension $mn$ (cf. \cite{Y2}).
I propose to name the space ${\mathcal P}_{n,m}$ the {\it Minkowski}-{\it Euclid space}
because H. Minkowski \cite{Mi} found a fundamental domain for $\mathcal P_n$ with respect to the
arithmetic subgroup $GL(n,\BZ)$ by means of the reduction theory. In this setting,
using the invariant differential operators on ${\mathcal P}_{n,m}$ we may develop the theory of automorphic
forms on ${\mathcal P}_{n,m}$ generalizing that on $\mathcal P_n$.

\vskip 0.2cm The aim of this paper is to study differential
operators on ${\mathcal P}_{n,m}$ which are invariant under the
action (1.2) of $\Gnm$. This article is organized as follows. In
Section 2, we review differential operators on $\CP$ invariant
under the action (1.1) of $GL(n,\BR)$. In Section 3, we
investigate differential operators on ${\mathcal P}_{n,m}$
invariant under the action (1.2) of $\Gnm$.
For two positive integers $m$ and $n$, we let
\begin{equation*}
S_{n,m}=\,\left\{\, (X,Z)\,|\ X=\,{}^tX\in \BR^{(n,n)},\ Z\in\BR^{(m,n)}\,\right\}
\end{equation*}
be the real vector space of dimension ${{n(n+1)}\over 2}+mn.$
From the adjoint action of the group $GL_{n,m}$,
we have the {\it natural action} of the orthogonal group $O(n)$
on $S_{n,m}$ given by
\begin{equation}
k\cdot (X,Z)=\,( k\,X\,{}^t k, Z\,{}^tk),\quad k\in O(n),\ (X,Z)\in S_{n,m}.
\end{equation}
The action (1.3) of $K=O(n)$ induces canonically the representation $\sigma$ of $O(n)$ on the polynomial algebra
$\textrm{Pol}(S_{n,m})$ consisting of complex valued polynomial functions on $S_{n,m}.$
Let $ \textrm{Pol}(S_{n,m})^K$ denote the subalgebra of $ \textrm{Pol}(S_{n,m})$ consisting of all polynomials on
$S_{n,m}$ invariant under the representation $\sigma$ of $O(n)$, and
$\BD ({\mathcal P}_{n,m})$ denote the algebra of all differential operators on ${\mathcal P}_{n,m}$ invariant under the
action $(1.2)$ of $GL_{n,m}$. We see that there is a canonically defined linear bijection of $\textrm{Pol}(S_{n,m})^K$
onto $\BD ({\mathcal P}_{n,m})$ which is not multiplicative. We will see that $\BD ({\mathcal P}_{n,m})$ is {\it not} commutative.
The main important problem is to find a complete list of explicit generators of $ \textrm{Pol}(S_{n,m})^K$
and a complete list of explicit generators of $\BD ({\mathcal P}_{n,m})$. We propose several natural problems.
We present some explicit invariant differential
operators which might be useful. In Section 4, we deal with the case $n=1$. In Section 5, we deal with the case $n=2$ and $m=1,2$.
In Section 6, we deal with the case $n=3$ and $m=1,2$. In Section 7, we deal with the case $n=4$ and $m=1,2$.
In the final section, we present some open problems and discuss a notion of automrphic forms on ${\mathcal P}_{n,m}$
using $GL_{n,m}$-invariant differential operators on the Minkowski-Euclid space ${\mathcal P}_{n,m}.$

\vskip 0.5cm
\noindent {\bf Acknowledgements:} This work was in part done during my stay at the Max-Planck-Institut f{\"u}r Mathematik in
Bonn. I am very grateful for the hospitality and financial support.
I would like to give my hearty and deep thanks to Minoru Itoh for his
interest in this work and many fruitful discussions.

\vskip 0.31cm \noindent {\bf Notations:} \ \ We denote by
$\BQ,\,\BR$ and $\BC$ the field of rational numbers, the field of
real numbers and the field of complex numbers respectively. We
denote by $\BZ$ and $\BZ^+$ the ring of integers and the set of
all positive integers respectively. The symbol ``:='' means that
the expression on the right is the definition of that on the left.
For two positive integers $k$ and $l$, $F^{(k,l)}$ denotes the set
of all $k\times l$ matrices with entries in a commutative ring
$F$. For a square matrix $A\in F^{(k,k)}$ of degree $k$,
$\textrm{tr}(A)$ denotes the trace of $A$. For any $M\in F^{(k,l)},\
^t\!M$ denotes the transposed matrix of $M$. $I_n$ denotes the
identity matrix of degree $n$. For a positive integer $n$, $I_n$ denotes the
identity matrix of degree $n$.

\end{section}


\begin{section}{{\bf Review on Invariant Differential Operators on $\CP$}}
\setcounter{equation}{0}

For a variable $Y=(y_{ij})\in \CP$, we set   $$dY=(dy_{ij})\quad
\textrm{and}\quad \Yd=\left( { {1+\delta_{ij}}\over 2 } {
{\partial}\over {\partial y_{ij} }}\right),$$ where $\delta_{ij}$
denotes the Kronecker delta symbol.

\vskip 0.2cm For a fixed element $g\in GL(n,\BR)$, we put
$$Y_*=g\cdot Y=gY\,^tg,\quad Y\in \CP.$$
Then
\begin{equation}
dY_*=g\,dY\,^tg \quad \textrm{and}\quad {{\partial}\over {\partial
Y_*}}=\,^tg^{-1} \Yd\, g^{-1}.
\end{equation}

We consider the following differential operators
\begin{equation}
D_i= \textrm{tr}\left( \left( Y\Yd \right)^i\right),\quad
i=1,2,\cdots,n,
\end{equation}

\noindent where $ \textrm{tr}(A)$ denotes the trace of a square
matrix $A$. By Formula (2.1), we get
\begin{equation*}
\left( Y_* {{\partial}\over {\partial Y_*}}\right)^i=g\,\left(
Y\Yd\right)^i g^{-1}
\end{equation*}

\noindent for any $g\in GL(n,\BR)$. So each $D_i$ is invariant
under the action (1.1) of $GL(n,\BR)$.

\vskip 0.2cm Selberg [9] proved the following.

\begin{theorem}
The algebra ${\mathbb D}(\CP)$ of all differential operators on
$\CP$ invariant under the action (1.1) of $GL(n,\BR)$ is generated
by $D_1,D_2,\cdots,D_n.$ Furthermore $D_1,D_2,\cdots,$ \par
\noindent $D_n$ are algebraically independent and ${\mathbb
D}(\CP)$ is isomorphic to the commutative ring \par\noindent $\BC
[x_1,x_2,\cdots,x_n]$ with $n$ indeterminates
$x_1,x_2,\cdots,x_n.$
\end{theorem}

\begin{proof} The proof can be found in [4], p.\,337, [7], pp.\,64-66 and  [10],
pp.\,29-30. The last statement follows immediately from the work
of Harish-Chandra [1,\,2] or [4], p.\,294. \end{proof}

Let $\fg=\BR^{(n,n)}$ be the Lie algebra of $GL(n,\BR)$. The
adjoint representation $ \textrm{Ad}$ of $GL(n,\BR)$ is given by
\begin{equation*}
\textrm{Ad}(g)=gXg^{-1},\quad g\in GL(n,\BR),\ X\in\fg.
\end{equation*}

\noindent The Killing form $B$ of $\fg$ is given by
\begin{equation*}
B(X,Y)= 2n\,\textrm{tr}(XY)-2\,\textrm{tr}(X)\,\textrm{tr}(Y),\quad X,Y\in \fg.
\end{equation*}

\noindent Since $B(aI_n,X)=0$ for all $a\in \BR$ and $X\in \frak g,\ B$ is degenerate. So
the Lie algebra ${\frak g}$ of $GL(n,\BR)$ is not semi-simple.

\vskip 0.3cm
The Lie algebra $\fk$ of $K$ is
\begin{equation*}
\fk=\left\{\,X\in \fg\,\vert\ X+\,^tX=0\ \right\}.\end{equation*}
We let $\fp$ be the subspace of $\fg$ defined by
\begin{equation*}
\fp=\left\{\,X\in \fg\,\vert\ X=\,^tX\in \BR^{(n,n)}\,
\right\}.\end{equation*}

\noindent Then
$$\fg=\fk\oplus \fp$$

\noindent is the direct sum of $\fk$ and $\fp$ with respect to the
Killing form $B$. Since $ \textrm{Ad}(k)\fp \subset \fp$ for any
$k\in K,\ K$ acts on $\fp$ via the adjoint representation by
\begin{equation}
k\cdot X= \textrm{Ad}(k)X=kX\,^tk,\quad k\in K,\ X\in \fp.
\end{equation}

The action (2.3) induces the action of $K$ on the polynomial
algebra $ \textrm{Pol}(\fp)$ of $\fp$ and the symmetric algebra
$S(\fp)$. We denote by $ \textrm{Pol}(\fp)^K$ (resp. $S(\fp)^K$)
the subalgebra of $ \textrm{Pol}(\fp)$ (resp. $S(\fp)$) consisting
of all $K$-invariants. The following inner product $(\ ,\ )$ on
$\fp$ defined by $$(X,Y)=B(X,Y),\quad X,Y\in \fp$$

\noindent gives an isomorphism as vector spaces
\begin{equation}
\fp\cong \fp^*,\quad X\mapsto f_X,\quad X\in \fp,
\end{equation}

\noindent where $\fp^*$ denotes the dual space of $\fp$ and $f_X$
is the linear functional on $\fp$ defined by
$$f_X(Y)=(Y,X),\quad Y\in \fp.$$

\noindent It is known that there is a canonical linear bijection
of $S(\fp)^K$ onto ${\mathbb D}(\CP)$. Identifying $\fp$ with
$\fp^*$ by the above isomorphism (2.4), we get a canonical linear
bijection
\begin{equation}
\Theta_n:\textrm{Pol}(\fp)^K \lrt {\mathbb D}(\CP)
\end{equation}

\noindent of $ \textrm{Pol}(\fp)^K$ onto ${\mathbb D}(\CP)$. The
map $\Theta_n$ is described explicitly as follows. We put
$N=n(n+1)/2$. Let $\left\{ \xi_{\alpha}\,|\ 1\leq \alpha \leq N\,
\right\}$ be a basis of $\fp$. If $P\in \textrm{Pol}(\fp)^K$, then
\begin{equation}
\Big(\Theta_n (P)f\Big)(gK)=\left[ P\left( {{\partial}\over {\partial
t_{\al}}}\right)f\left(g\,\text{exp}\, \left(\sum_{\al=1}^N
t_{\al}\xi_{\al}\right) K\right)\right]_{(t_{\al})=0},
\end{equation} where $f\in C^{\infty}({\mathcal P}_{n})$. We refer to [3,\,4] for more detail. In
general, it is very hard to express $\Theta_n (P)$ explicitly for a
polynomial $P\in \textrm{Pol}(\fp)^K$.

\skip 0.52cm We let
\begin{equation}
q_i(X)= \textrm{tr}\big(X^i\big),\quad i=1,2,\cdots,n
\end{equation}

\noindent be the polynomials on $\fp$. Here we take a coordinate
$x_{11}, x_{12},\cdots,x_{nn}$ in $\fp$ given by

$$X=\begin{pmatrix}
  x_{11} & {\frac 12}x_{12}  & \hdots &
{\frac 12}x_{1n} \\ {\frac 12}x_{12} & x_{22} & \hdots & {\frac
12}x_{2n}\\
\vdots& \vdots & \ddots&\vdots\\
{\frac 12}x_{1n} & {\frac 12}x_{2n} &\hdots &x_{nn}
\end{pmatrix}.$$

\noindent For any $k\in K$,
\begin{equation*}
(k\cdot q_i)(X)=q_i(k^{-1}Xk)=\textrm{tr}(k^{-1}X^ik)=q_i(X),\quad
i=1,2,\cdots,n.
\end{equation*}

\noindent Thus $q_i\in \textrm{Pol}(\fp)^K$ for $ i=1,2,\cdots,
n.$ By a classical invariant theory (cf. [5,\,11]), we can prove
that the algebra $\textrm{Pol}(\fp)^K$ is generated by the
polynomials $q_1,q_2,\cdots, q_n$ and that $q_1,q_2,\cdots, q_n$
are algebraically independent. Using Formula (2.6), we can show
without difficulty that

\begin{equation*}
\Theta_n(q_1)= \textrm{tr} \left( 2 Y \Yd \right).
\end{equation*}

\noindent However $\Theta_n(q_i)\ (i=2,3,\cdots,n)$ are still not
known explicitly.

\vskip 0.2cm We propose the following conjecture.

\vskip 0.2cm \noindent
$ \textbf{Conjecture 1.}$ $ \textit{For\ any
}\ n,$
\begin{equation*}
\Theta_n(q_i)= \textrm{tr} \left( \left(2 Y \Yd\right)^i \right),\quad
i=1,2,\cdots,n.
\end{equation*}

\vskip 0.3cm \noindent
$ \textbf{Remark.}$ The author checked that
the above conjecture is true for $n=1,2.$

\vskip 0.2cm For a positive real number $A$,
\begin{equation*}
ds_{n;A}^2= A\cdot \textrm{tr}\big(Y^{-1}dY\,Y^{-1} dY\big)
\end{equation*}

\noindent is a Riemannian metric on $\CP$ invariant under the
action (1.1). The Laplacian $\Delta_{n;A}$ of $ds_{n;A}^2$ is
\begin{equation*}
\Delta_{n;A}={\frac 1A}\, \textrm{tr} \left( \left( Y \Yd
\right)^2\right).
\end{equation*}

\newcommand\da{ {\partial}\over {\partial y_1} }
\newcommand\db{ {\partial}\over {\partial y_2} }
\newcommand\dc{ {\partial}\over {\partial y_3} }
For instance, we consider the case $n=2$ and $A>0$. If we write
for $Y\in {\mathcal P}_2$
$$Y=\begin{pmatrix} y_1 & y_3 \\ y_3 & y_2 \end{pmatrix}\quad
\text{and}\quad {{\partial}\over {\partial Y}}=\begin{pmatrix} \da & {\frac 12} {\dc} \\
{\frac 12} {\dc} & \db \end{pmatrix},$$ then

\begin{eqnarray*}
 ds^2_{2;A}&=& A\,\, \textrm{tr}\big(Y^{-1}dY\,Y^{-1}dY\big)\cr
 &=&{\frac A{\big(y_1y_2-y^2_3\big)^2}}\,\,  \Big\{ y_2^2\,
 dy_1^2\,+\,y_1^2\,dy_2^2\,+\,2\,\big(y_1y_2+y_3^2\big)\,dy_3^2 \cr & & \ \ \ \
 \ \ \
 +\,2\,y_3^2 \,\,dy_1dy_2- 4\,y_2\,y_3\,\,dy_1dy_3 -4\,y_1y_3\,\,dy_2 dy_3
\Big\}
\end{eqnarray*}

and its Laplacian $\Delta_{2;A}$ on ${\mathcal P}_2$ is

\newcommand\POB{ {{\partial}\over
{\partial{\overline \Omega}}} }
\newcommand\PZB{ {{\partial}\over {\partial{\overline Z}}} }
\newcommand\PX{ {{\partial}\over{\partial X}} }
\newcommand\PY{ {{\partial}\over {\partial Y}} }
\newcommand\PU{ {{\partial}\over{\partial U}} }
\newcommand\PV{ {{\partial}\over{\partial V}} }
\newcommand\PO{ {{\partial}\over{\partial \Omega}} }
\newcommand\PZ{ {{\partial}\over{\partial Z}} }
\begin{eqnarray*}
 \Delta_{2;A}&=& {\frac 1A}\,\, \textrm{tr}  \left( \left( Y\PY\right)^2\right)\cr
&=&{\frac 1A}\, \Biggl\{ y_1^2\,{{\partial^2}\over{\partial
y_1^2}}+y_2^2\,{{\partial^2}\over{\partial y_2^2}}+ {\frac
12}(y_1y_2+y_3^2){{\partial^2}\over{\partial y_3^2}}\cr & &\ \ \
+2\left( y_3^2 {{\partial^2}\over {\partial y_1\partial y_2}}+
y_1y_3 {{\partial^2}\over {\partial y_1\partial y_3}}+y_2y_3
{{\partial^2}\over {\partial y_2\partial y_3}} \right)\cr & &\ \ \
+ {3 \over 2}\left( y_1{ {\partial}\over {\partial y_1} }+y_2{
{\partial}\over {\partial y_2} }+y_3{ {\partial}\over {\partial
y_3} }\right)\Biggl\}.
\end{eqnarray*}

\end{section}

\vskip 1cm


\begin{section}{{\bf Invariant Differential Operators on ${\mathcal P}_{n,m}$}}
\setcounter{equation}{0}

\vskip 0.3cm For a variable $(Y,V)\in {\mathcal P}_{n,m}$ with
$Y\in \CP$ and $V\in \Rmn$, we put
$$Y=(y_{ij})\ \text{with}\ y_{ij}=y_{ji},\ \
V=(v_{kl}),$$   $$ dY=(dy_{ij}),\ \ dV=(dv_{kl}),$$
$$[dY]=\wedge_{i\leq j}dy_{ij},\ \ \ \ \
[dV]=\wedge_{k,l}dv_{kl},$$ and
$$\Yd=\left( { {1+\delta_{ij}}\over 2 } {
{\partial}\over {\partial y_{ij} }}\right),\ \ \ \Vd=\left({
{\partial}\over {\partial v_{kl}} } \right),$$ where $1\leq
i,j,l\leq n$ and $1\leq k\leq m.$

\vskip 0.2cm For a fixed element $(g,\lambda)\in \Gnm,$ we write
\begin{equation*}
(Y_{\star},V_{\star})=(g,\lambda)\cdot
(Y,V)=\big(g\,Y\,^tg,(V+\lambda)\,^tg\big),
\end{equation*}

\noindent where $(Y,V)\in {\mathcal P}_{n,m}$. Then we get

\begin{equation}
Y_{\star}= g\,Y\,\,^tg,\quad V_{\star}=(V+\lambda)\,^tg
\end{equation}
\noindent and

\begin{equation}
{ {\partial }\over {\partial Y_{\star}} }=
\,^tg^{-1}\,\Yd\,g^{-1},\quad {{\partial}\over {\partial
V_{\star}} } ={{\partial}\over {\partial V}}\,g^{-1}.
\end{equation}

\newcommand\fgt{\fg_{\star}}
\newcommand\Adt{$ \textrm{Ad}_{\star}$ }
\newcommand\fpt{\fp_{\star}}
\newcommand\Yst{Y_{\star}}
\newcommand\Vst{V_{\star}}
\newcommand\dYt{ { {\partial}\over {\partial \Yst} } }
\newcommand\dVt{ { {\partial}\over {\partial \Vst} } }

\begin{lemma}
For all two positive real numbers $A$ and $B$,
the following metric $ds^2_{n,m;A,B}$ on ${\mathcal P}_{n,m}$ defined by
\begin{equation}
ds^2_{n,m;A,B}=\,A\,\s(Y^{-1}dY\,Y^{-1}dY)\,+\,B\,\s(Y^{-1}\,^t(dV)\,dV)
\end{equation}
is a Riemannian metric on ${\mathcal P}_{n,m}$ which is invariant under the action (1.2) of $\Gnm$. The
Laplacian $\Delta_{n,m;A,B}$ of $({\mathcal P}_{n,m},\,ds^2_{n,m;A,B})$ is given by
$$\Delta_{n,m;A,B}=\,{\frac 1A}\,\s \left(\left( Y {{\partial}\over {\partial
Y}}\right)^2\right) -{\frac m{2\,A}}\,\s\left( Y {{\partial}\over
{\partial Y}}\right)\,+ \,{\frac 1B}\,\sum_{k\leq p} \left( \left({{\partial}\over
{\partial V}}\right) Y
{}^{{}^{{}^{{}^\text{\scriptsize $t$}}}}\!\!\!
\left({{\partial}\over {\partial
V}}\right)\right)_{kp}.$$ Moreover $\Delta_{n,m;A,B}$ is a
differential operator of order 2 which is invariant under the
action (1.2) of $\Gnm.$
\end{lemma}

\vskip 0.3cm\noindent
{\it Proof.} The proof can be found in \cite{Y2}.
\hfill $\square$

\vskip 0.5cm
\begin{lemma} The following volume element $dv_{n,m}(Y,V)$ on ${\mathcal P}_{n,m}$ defined by
\begin{equation}
dv_{n,m}(Y,V)=(\det\,Y)^{-{ {n+m+1}\over 2} }[dY][dV]
\end{equation} is invariant under the action (1.2)
of $\Gnm$.
\end{lemma}

\vskip 0.3cm\noindent
{\it Proof.} The proof can be found in \cite{Y2}.
\hfill $\square$

\vskip 0.5cm

\begin{theorem} Any geodesic through the origin $(I_n,0)$ for the Riemannian metric $ds^2_{n,m;1,1}$
is of the form
$$\gamma(t)=\left(\,\la(2t)[k],\,Z\left(\int_0^t
\la (t-s)ds\right)[k]\right),$$ where $k$ is a fixed element of
$O(n),\ Z$ is a fixed $h\times g$ real matrix, $t$ is a real
variable, $\la_1,\la_2,\cdots,\la_n$ are fixed real numbers but
not all zero and
$$\la(t):=\text{diag}\,(e^{\la_1t},\cdots,e^{\la_n t}).$$
Furthermore, the tangent vector $\g'(0)$ of the geodesic $\g(t)$
at $(I_n,0)$ is $(D[k],Z),$ where
$D=\text{diag}\,(2\la_1,\cdots,2\la_n).$
\end{theorem}

\vskip 0.3cm\noindent
{\it Proof.} The proof can be found in \cite{Y2}.
\hfill $\square$

\vskip 0.5cm
\begin{theorem}
Let $(Y_0,V_0)$ and $(Y_1,V_1)$ be two points in ${\mathcal P}_{n,m}.$ Let $g$ be an element in $GL(n,\BR)$
such that $Y_0[\,{}^tg]=\,I_n$ and $Y_1[\,{}^tg]$ is diagonal.
Then the length $s\big((Y_0,V_0),(Y_1,V_1) \big)$ of the geodesic
joining $(Y_0,V_0)$ and $(Y_1,V_1)$ for the $GL_{n,m}$-invariant Riemannian metric
$ds_{n,m;A,B}^2$ is given by
\begin{equation}
s\big((Y_0,V_0),(Y_1,V_1)\big)=\,A\,\left\{ \sum_{j=1}^n (\ln t_j)^2 \right\}^{1/2}\,+\,B\,
\int_0^1 \left( \sum_{j=1}^n \Delta_j\,e^{-(\ln t_j)\,t}\right)^{1/2}\,dt,
\end{equation}
where $\Delta_j=\,\sum_{k=1}^m{\widetilde v}_{kj}^2 \ (1\leq j\leq n)$ with
$(V_1-V_0)\,{}^tg=\,({\widetilde v}_{kj})$ and
$t_1,\cdots,t_n$ denotes the zeros of $\det (t\,Y_0-Y_1)$.
\end{theorem}

\vskip 0.3cm\noindent
{\it Proof.} The proof can be found in \cite{Y2}.
\hfill $\square$

\vskip 0.53cm
The Lie algebra $\fg_{\star}$ of $\Gnm$ is given by
$$\fg_{\star}=\Big\{\,(X,Z)\,\vert\ X\in \BR^{(n,n)},\ Z\in
\BR^{(m,n)}\,\Big\} $$ equipped with the following Lie bracket
\begin{equation*}
\big[(X_1,Z_1),\,(X_2,Z_2)\big]=\big([X_1,X_2]_0,\,Z_2\,^tX_1-Z_1\,^tX_2\big),\end{equation*}
where $[X_1,X_2]_0=X_1X_2-X_2X_1$ denotes the usual matrix bracket
and $(X_1,Z_1),\,(X_2,Z_2)\in \fgt$. The adjoint representation
$\text{Ad}_{\star}$ of $\Gnm$ is given by
\begin{equation}
\text{Ad}_{\star}\,((g,\lambda))(X,Z)=\big(\,
gXg^{-1},\,(Z-\lambda\,^tX)\,^tg\,\big),\end{equation} where
$(g,\lambda)\in \Gnm$ and $(X,Z)\in \fgt$. And the adjoint
representation $\text{ad}_{\star}$ of $\fgt$ on
$\text{End}\,(\fgt)$ is given by
$$\text{ad}_{\star}\,\big((X,Z)\big)\big((X_1,Z_1)\big)=\big[(X,Z),(X_1,Z_1)\big].$$ We see
that the Killing form $B_{\star}$ of $\fgt$ is given by
$$B_{\star}\big((X_1,Z_1),(X_2,Z_2)\big)=(2\,n+m)\, \textrm{tr} (X_1X_2)-2\, \textrm{tr}(X_1)\,\textrm{tr}(X_2).$$ The
Lie algebra $\fk$ of $K$ is
\begin{equation*}
\fk=\Big\{\,(X,0)\in \fgt\,\big\vert\ X+\,^tX=0\
\Big\}.\end{equation*} We let $\fpt$ be the subspace of $\fgt$
defined by
\begin{equation*}
\fpt=\Big\{\,(X,Z)\in \fgt\,\big\vert\ X=\,^tX\in \BR^{(n,n)},\
Z\in \Rmn\, \Big\}.\end{equation*} Then we have the following
relation
\begin{equation*}
[\fk,\fk]\subset \fk\ \ \ \ \ \ \text{and}\ \ \ \ \ \
[\fk,\fpt]\subset \fpt. \end{equation*} In addition, we have
$$\fgt=\fk\oplus \fpt\ \ \ \ \ \text{(\,the\ direct\ sum\,)}.$$

$K$ acts on $\fpt$ via the adjoint representation $
\textrm{Ad}_{\star}$ of $\Gnm$ by

\begin{equation}
k\cdot (X,Z)=\big(kX\,^tk, Z\,^tk\big),\quad k\in K,\ (X,Z)\in
\fpt.
\end{equation}

\noindent The action (3.7) induces the action of $K$ on the
polynomial algebra $ \textrm{Pol}(\fpt)$ of $\fpt$ and the
symmetric algebra $S(\fpt)$. We denote by $ \textrm{Pol}(\fpt)^K$
(resp. $S(\fpt)^K$) the subalgebra of $ \textrm{Pol}(\fpt)$ (resp.
$S(\fpt)$) consisting of all $K$-invariants. The following inner
product $(\ ,\ )_{\star}$ on $\fpt$ defined by
$$\big( (X_1,Z_1),(X_2,Z_2)\big)_{\star}= \textrm{tr}(X_1X_2)+ \textrm{tr}(Z_1\,^tZ_2),
\quad (X_1,Z_1),(X_2,Y_2)\in \fpt$$

\noindent gives an isomorphism as vector spaces
\begin{equation}
\fpt\cong \fpt^*,\quad (X,Z)\mapsto f_{X,Z},\quad (X,Z)\in \fpt,
\end{equation}

\noindent where $\fpt^*$ denotes the dual space of $\fpt$ and
$f_{X,Z}$ is the linear functional on $\fpt$ defined by
$$f_{X,Z}\big( (X_1,Z_1)\big)=\big( (X,Z),(X_1,Z_1)\big)_{\star},\quad (X_1,Z_1)\in \fpt.$$

\noindent Let ${\mathbb D}({\mathcal P}_{n,m})$ be the algebra of
all differential operators on ${\mathcal P}_{n,m}$ that are
invariant under the action (1.2) of $\Gnm$. It is known that there
is a canonical linear bijection of $S(\fpt)^K$ onto ${\mathbb
D}({\mathcal P}_{n,m})$. Identifying $\fpt$ with $\fpt^*$ by the
above isomorphism (3.5), we get a canonical linear bijection
\begin{equation}
\Theta_{n,m}:\textrm{Pol}(\fpt)^K \lrt {\mathbb D}({\mathcal P}_{n,m})
\end{equation}

\noindent of $ \textrm{Pol}(\fpt)^K$ onto ${\mathbb D}({\mathcal
P}_{n,m})$. The map $\Theta_{n,m}$ is described explicitly as follows.
We put $N_{\star}=n(n+1)/2+mn$. Let $\left\{ \eta_{\alpha}\,|\
1\leq \alpha \leq N_{\star}\, \right\}$ be a basis of $\fpt$. If
$P\in \textrm{Pol}(\fpt)^K$, then
\begin{equation}
\Big(\Theta_{n,m} (P)f\Big)(gK)=\left[ P\left( {{\partial}\over
{\partial t_{\al}}}\right)f\left(g\,\text{exp}\,
\left(\sum_{\al=1}^{N_{\star}} t_{\al}\eta_{\al}\right)
K\right)\right]_{(t_{\al})=0}, \end{equation} where $f\in
C^{\infty}({\mathcal P}_{n,m})$. We refer to [4], pp.\,280-289. In
general, it is very hard to express $\Theta_{n,m}(P)$ explicitly for a
polynomial $P\in \textrm{Pol}(\fpt)^K$.

\vskip 0.2cm
We take a coordinate $(X,Z)$ in $\fpt$ such
that
$$X=\begin{pmatrix}
  x_{11} & {\frac 12}x_{12}  & \hdots &
{\frac 12}x_{1n} \\ {\frac 12}x_{12} & x_{22} & \hdots & {\frac
12}x_{2n}\\
\vdots& \vdots & \ddots&\vdots\\
{\frac 12}x_{1n} & {\frac 12}x_{2n} &\hdots &x_{nn}
\end{pmatrix}\in \fp \quad \textrm{and}\quad Z=(z_{kl})\in \BR^{(m,n)}.$$

We define the polynomials $\alpha_j,\ \beta_{pq}^{(k)},\ R_{jp}$ and $S_{jp}$ on
$\fpt$ by

\begin{eqnarray}
\alpha_j (X,Z)&=& \textrm{tr}   \big( X^j\,\big),\quad 1\leq j\leq n,\\
\beta_{pq}^{(k)}(X,Z)&=&  \big(Z\,X^k\,^tZ\,\big)_{pq},\quad 0\leq k\leq n-1,\ 1\leq
p\leq q\leq m,\\
R_{jp}(X,Z)&=& \textrm{tr}\big (X^j(\,^tZZ\,)^p\,\big),\quad 0\leq
j\leq n-1,\ 1\leq p\leq m,\\
S_{jp}(X,Z)&=& \det \big (X^j(\,^tZZ\,)^p\,\big),\quad 0\leq
j\leq n-1,\ 1\leq p\leq m,
\end{eqnarray}

\noindent where $\big(Z\,^t\!Z\big)_{pq}$ $\big( \textrm{resp.}\
\big( ZX\,^tZ\big)_{pq}\big)$ denotes the $(p,q)$-entry of
$Z\,^t\!Z$ $\big( \textrm{resp.}\ ZX\,^tZ\big)$.

\vskip 0.2cm For any $m\times m$ real matrix $S$, we define the
polynomials $M_{j;S},\ Q_{p;S}$ and $R_{i,p,j;S}$ on $\fpt$ by

\begin{eqnarray}
M_{j;S} (X,Z)&=& \textrm{tr}   \big( (X+\,^tZSZ)^j\,\big),\quad 1\leq j\leq n,\\
Q_{p;S}(X,Z)&=& \textrm{tr} \big( (\,^tZ\,S\,Z)^p\,\big),\quad
1\leq p\leq n,\\
\Omega_{i,p,j;S}(X,Z)&=&\, \textrm{tr}\Big( X^i (\,^tZSZ)^p
(X+\,^tZSZ)^j\,\Big),\\
\Theta_{i,p,j;S}(X,Z)&=&\, \det \Big( X^i (\,^tZSZ)^p
(X+\,^tZSZ)^j\,\Big),
\end{eqnarray}

\noindent where $0\leq i,j\leq n-1,\ 1\leq p\leq n.$ We see that all
$\alpha_j,\ \beta_{pq}^{(k)},\ R_{jp},\ S_{jp},\  M_{j;S},\ Q_{p;S},\ \Omega_{i,p,j;S} $
and $\Theta_{i,p,j;S}$ are elements of $\textrm{Pol}(\fpt)^K$.

 \vskip 0.355cm We propose the following natural problems.

\vskip 0.2cm \noindent $ \textbf{Problem 1.}$ Find a complete list of explicit generators
of $\textrm{Pol}(\fpt)^K$.

\vskip 0.32cm \noindent $ \textbf{Problem 2.}$ Find all the relations among a set of generators of $\textrm{Pol}(\fpt)^K$.

\vskip 0.32cm \noindent $ \textbf{Problem 3.}$ Find an easy or effective way to
express the images of the above invariant polynomials under the
Helgason map $\Theta_{n,m}$ explicitly.

\vskip 0.32cm \noindent $ \textbf{Problem 4.}$ Decompose $\textrm{Pol}(\fpt)^K$ into
$O(n)$-irreducibles.

\vskip 0.32cm \noindent $ \textbf{Problem 5.}$ Find a complete list of explicit generators of the algebra $\BD({\mathcal P}_{n,m})$.
Or construct explicit $GL_{n,m}$-invariant differential operators on ${\mathcal P}_{n,m}.$

\vskip 0.32cm \noindent $ \textbf{Problem 6.}$ Find all the relations among a set of generators of $\BD({\mathcal P}_{n,m})$.

\vskip 0.32cm \noindent $ \textbf{Problem 7.}$ Is $\textrm{Pol}(\fpt)^K$ finitely generated ? Is $\BD({\mathcal P}_{n,m})$
finitely generated ?

\vskip 0.5cm
\begin{theorem} $\textrm{Pol}(\fpt)^K$ is generated by $\alpha_j\ (1 \leq j \leq n)$ and $\beta_{pq}^{(k)}$
($0\leq k\leq n-1,\ 1\leq p\leq q\leq m).$
\end{theorem}

\vskip 0.53cm We present some invariant differential operators on
${\mathcal P}_{n,m}.$ We define the differential operators $D_j,\
 \Omega_{pq}$ and $L_p$ on ${\mathcal P}_{n,m}$ by
\begin{equation}
D_j=\,\textrm{tr}\left(\left(\,2Y\Yd\right)^j\right),\ \ \ \ 1\leq
j\leq n,\end{equation}
\begin{equation}
\Omega_{pq}^{(k)}=\,\left\{\,{{\partial}\over{\partial V}}\left(
2Y{{\partial}\over{\partial Y}}\right)^kY\,{}^{{}^{{}^{{}^\text{\scriptsize $t$}}}}\!\!\!\left({{\partial}\over
{\partial V}}\right)\,\right\}_{pq},\quad 0\leq k\leq n-1,\ 1\leq p\leq q\leq
m\end{equation} and
\begin{equation}
L_p=\,\textrm{tr}\left(\left\{Y\,{}^{{}^{{}^{{}^\text{\scriptsize $t$}}}}\!\!\!\left(\Vd\right)\Vd\right\}^p\right),\
\ \ 1\leq p\leq m.\end{equation}
\noindent
Here for a matrix $A$ we denote by $A_{pq}$ the $(p,q)$-entry of $A$.

\vskip 0.3cm
Also we define the differential
operators $S_{jp}$ by
\begin{equation}
S_{jp}=\,\textrm{tr}\left(\left(\,2Y\Yd\right)^j\left\{Y\,{}^{{}^{{}^{{}^\text{\scriptsize $t$}}}}\!\!\!\left(\Vd\right)
\Vd\right\}^p\right),\end{equation} where $1\leq j\leq n$ and
$1\leq p\leq m.$ \vskip 0.12cm For any real matrix $S$ of degree
$m$, we define the differential operators $\Phi_{j;S},\ L_{p;S}$
and $\Phi_{i,p,j;S}$ by
\begin{equation}
\Phi_{j;S}=  \textrm{tr}
\left(\,\left\{\,Y\left(\,2\,\Yd+\,{}^{{}^{{}^{{}^\text{\scriptsize $t$}}}}\!\!\!\left(\Vd\right)S
\left(\Vd\right)\right)\right\}^j\right),\ \ \ 1\leq j\leq
n,\end{equation}
\begin{equation}
L_{p;S}=\,\textrm{tr}\left(\,\left\{\,Y\,{}^{{}^{{}^{{}^\text{\scriptsize $t$}}}}\!\!\!\left(\Vd\right)S\left(\Vd\right)
\right\}^p\,\right),\ \ \ \ \ 1\leq p\leq m\end{equation} and
\begin{equation}
\Phi_{i,p,j;S}(X,Z) =
\end{equation}
\begin{equation*}
\,\textrm{tr}\left(\,\left(\,2\,Y\Yd\right)^i \left(
Y\,{}^{{}^{{}^{{}^\text{\scriptsize $t$}}}}\!\!\!\left(\Vd\right) S\left(\Vd\right) \right)^p \left\{
Y\left(\,2\,\Yd+\,{}^{{}^{{}^{{}^\text{\scriptsize $t$}}}}\!\!\!\left(\Vd\right)S\left(\Vd\right)
\right)\right\}^j\right). \end{equation*}

\vskip 0.5cm We want to mention the special invariant differential
operator on ${\mathcal P}_{n,m}$. In \cite{Y1}, the author studied the
following differential operator $M_{n,m,{\mathcal M}}$ on
${\mathcal P}_{n,m}$ defined by
\begin{equation}
M_{n,m,\CM}=\text{det}\,(Y)\cdot\text{det}\,\left( \Yd+{1\over
{8\pi}}\, {}^{{}^{{}^{{}^\text{\scriptsize $t$}}}}\!\!\!\left(\Vd\right)\CM^{-1}\left(\Vd\right)\right),
\end{equation} where $\CM$ is a positive definite, symmetric
half-integral matrix of degree $m$. This differential operator
characterizes {\it singular\ Jacobi\ forms}. For more detail, we
refer to \cite{Y1}. According to (3.1) and (3.2), we see easily that
the differential operator $M_{n,m,\CM}$ is invariant under the
action (1.2) of $\Gnm$.

\vskip 0.3cm
\noindent
{\bf Question\,:} Calculate the inverse
of $M_{n,m,\CM}$ under the Helgason map $\Theta_{n,m}$.

\end{section}

\vskip 1cm


\begin{section}{{\bf The Case $n=1$}}
\setcounter{equation}{0}

\vskip 0.3cm In this section, we consider the case $n=m=1$ and the case $n=1$ and $m\geq 2$ separately.
\vskip 0.3cm\noindent
{\bf 4.1.} {\bf The Case $n=1$ and $m=1$}
\vskip 0.3cm In this case,
\begin{equation*}
GL_{1,1}=\BR^{\times}\ltimes \BR,\quad K=O(1), \quad
{\mathcal P}_{1,1}=\BR^+\times \BR,
\end{equation*}
\noindent where $\BR^{\times}=\{ \,a\in \BR\,|\ a\neq 0\,\}$ and
$\BR^+=\{\ a\in \BR\,|\ a>0\,\}.$ Clearly $\fk=0$ and
$\fp_{\star}=\fg_{\star}=\left\{ \,(x,z)\,|\ x,z\in\BR\, \right\}$.
Then $e=(1,0)$ and $f=(0,1)$ form the standard basis for $\fp_{\star}$. Using this
basis, we take
a coordinate $(x,z)$ in $\fp_{\star}$, that is,
if $w \in \fp_{\star}$, we write $w=xe +zf.$ We can show that $ \textrm{Pol}(\fp_{\star})^K$ is
generated by the following polynomials
\begin{equation*}
\alpha(x,z)=x\quad \textrm{and}\quad \beta(x,z)=z^2.
\end{equation*}
The generators $\alpha$ and $\beta$ are {\it algebraically independent}.
Let $(y,v)$ be a
coordinate in ${\mathcal P}_{1,1}$ with $y>0$ and $v\in\BR.$
Then using Formula
(3.10), we can show that
\begin{equation*}
\Theta_{1,1}(\alpha)=\,2y\,{{\partial}\over {\partial y}}\quad
\textrm{and}\quad \Theta_{1,1} (\beta)=\,y{{\partial^2}\over {\partial
v^2}}.
\end{equation*}
\noindent We see that $\Theta_{1,1}(\alpha)$ and $\Theta_{1,1}(\beta)$ generate the algebra ${\mathbb D}({\mathcal
P}_{1,1})$ and are {\it algebraically dependent}. Indeed we have the following noncommutation relation
$$\Theta_{1,1}(\alpha)\Theta_{1,1}(\beta)\,-\,\Theta_{1,1}(\beta)\Theta_{1,1}(\alpha)\,=\,2\,\Theta_{1,1}(\beta).$$ Hence the
algebra ${\mathbb D}({\mathcal P}_{1,1})$ is {\it not}
commutative. The unitary dual $\widehat K$ of $K$ consists of two elements. Let
\begin{equation*}
\textrm{Pol}(\fp_{\star})=\sum_{\tau\in{\widehat K}}m_\tau \tau
\end{equation*}
be the decomposition of $\textrm{Pol}(\fp_{\star})$ into $K$-irreducibles. It is easy to see that the multiplicity
$m_\tau$ of $\tau$ is infinite for all $\tau\in {\widehat K}.$ So the action of $K$ on $\textrm{Pol}(\fp_{\star})$
is not multiplicity-free. In this case the seven problems proposed in Section 3 are completely solved.


\vskip 0.5cm
\noindent
{\bf 4.2.} {\bf The Case $n=1$ and $m\geq 2$}
\vskip 0.3cm

We consider the case $n=1$ and $m\geq 2$. In this case,
\begin{equation*}
GL_{1,m}=\BR^{\times}\ltimes \BR^{(m,1)},\quad K=O(1), \quad
{\mathcal P}_{1,m}=\BR^+\times \BR^{(m,1)},\end{equation*}

\noindent where $\BR^{\times}=\{ \,a\in \BR\,|\ a\neq 0\,\}$ and
$\BR^+=\{\ a\in \BR\,|\ a>0\,\}.$ Clearly $\fk=0$ and
$\fp_{\star}=\fg_{\star}=\left\{ \,(x,z)\,|\ x\in\BR,\ z\in
\BR^{(m,1)}\ \right\}$. Let $\left\{ e_1,\cdots,e_m\right\}$ be
the standard basis of $\BR^{(m,1)}.$ Then
$$\eta_0=(1,0),\ \eta_1=(0,e_1),\ \eta_2=(0,e_2),\cdots,\ \eta_m=(0,e_m)$$
\noindent form a basis of $\fp_{\star}$. Using this basis, we take
a coordinate $(x,z_1,z_2,\cdots,z_m)$ in $\fp_{\star}$, that is,
if $w \in \fp_{\star}$, we write $w=x\eta_0 +\sum_{k=1}^m
z_k\eta_k.$ We can show that $ \textrm{Pol}(\fp_{\star})^K$ is
generated by the following polynomials
\begin{equation*}
\alpha(x,z)=x\quad \textrm{and}\quad \beta_{kl}(x,z)=z_k\, z_l,\quad 1\leq k\leq l\leq m,
\end{equation*}
\noindent where $z=(z_1,z_2,\cdots,z_m)$.
We see easily that one has the following relations
\begin{equation*}
\beta_{kk}\beta_{ll}=\beta_{kl}^2\qquad \textrm{for}\ 1\leq k < l\leq m
\end{equation*}
and
\begin{equation*}
\beta_{kk}\beta_{ll}^2\beta_{pp}=\beta_{kl}^2\beta_{lp}^2
\qquad \textrm{for}\ 1\leq k< l < p\leq m.
\end{equation*}
Therefore the generators $\alpha$ and $\beta_{kl}\,(1\leq k\leq l\leq m)$ are
{\it algebraically dependent}.

\vskip 0.2cm
Let $(y,v)$ be a
coordinate in ${\mathcal P}_{1,m}$ with $y>0$ and
$v=\,{}^t(v_1,v_2,\cdots,v_m)\in\BR^{(m,1)}.$ Then using Formula
(3.10), we can show that
\begin{equation*}
\Theta_{1,m}(\alpha)=\,2y\,{{\partial}\over {\partial y}}\quad
\textrm{and}\quad \Theta_{1,m} (\beta_{kl})=\,y{{\partial^2}\over {\partial
v_k\partial v_l}},\ \ \ 1\leq k\leq l\leq m.
\end{equation*}

\noindent We see that $\Theta_{1,m}(\alpha)$ and $\Theta_{1,m}(\beta_{kl})\,(1\leq
k\leq l\leq m)$ generate the algebra ${\mathbb D}({\mathcal
P}_{1,m})$. Although $\Theta_{1,m}(\beta_{kl})\,(1\leq k\leq l\leq m)$
commute with each other, $\Theta_{1,m}(\alpha)$ does not commute with any
$\Theta_{1,m}(\beta_{kl})$. Indeed, we have the noncommutation relation
$$\Theta_{1,m}(\alpha)\Theta_{1,m}(\beta_{kl})\,-\,\Theta_{1,m}(\beta_{kl})\Theta_{1,m}(\alpha)\,
=\,2\,\Theta_{1,m}(\beta_{kl}).$$ Hence the
algebra ${\mathbb D}({\mathcal P}_{1,m})$ is {\it not}
commutative. It is easily seen that the action of $K$ on $ \textrm{Pol}(\fp_{\star})$ is {\it not}
multiplicity-free.

\end{section}

\vskip 0.5cm


\begin{section}{{\bf The Case $n=2$}}
\setcounter{equation}{0}
\vskip 0.3cm In this section we deal with the case $n=2,\ m=1$ and the case $n=m=2.$

\vskip 0.5cm
\noindent
{\bf 5.1.} {\bf The Case $n=2$ and $m=1$}
\vskip 0.3cm
In this case, $$GL_{2,1}=GL(2,\BR)\ltimes
\BR^{(1,2)},\ \ K=O(2)\ \ \text{and} \ \ GL_{2,1}/K={\mathcal
P}_2\times \BR^{(1,2)}={\mathcal P}_{2,1}.$$ We see easily that
$$\fp_{\star}=\,\left\{\,(X,Z)\ \vert\ X=\,^tX\in \BR^{(2,2)},\ \ Z\in
\BR^{(1,2)} \,\right\}.$$ We put
$$e_1=\left(\begin{pmatrix} 1 & 0\\ 0 & 0 \end{pmatrix},0\right),\ \ \ e_{2}=\left(\begin{pmatrix} 0 &
0\\ 0 & 1 \end{pmatrix},0\right),\ \ \ e_3=\left(\begin{pmatrix} 0 & 1 \\
1 & 0
\end{pmatrix},0\right)$$ and $$f_1=(0,(1,0)),\ \ \ \ f_2=(0,(0,1)).$$ Then
$\{\,e_1,\,e_2,\,e_3,\,f_1,\,f_2\,\}$ forms a basis for
$\fp_{\star}.$
We write for variables $(X,Z)\in \fp_{\star}$ by $$X=\begin{pmatrix} x_1 & {\frac 12}x_3 \\
{\frac 12}x_3 & x_2 \end{pmatrix}\ \ \ \text{and} \ \ \
Z=(z_1,z_2).$$ The following polynomials
$$\alpha_1(X,Z)=\, \textrm{tr}   (X)=x_1\,+\,x_2,\ \ \ \ \ \alpha_2(X,Z)= \textrm{tr} (X^2)=
x_1^2\,+\,x_2^2\,+\,{\frac 12}\,x_3^2,$$
$$\xi(X,Z)=\,Z\,^tZ\,=\,z_1^2\,+\,z_2^2$$ and
$$\varphi(X,Z)=\,ZX\,^tZ\,=\,x_1\,z_1^2\,+\,x_2\,z_2^2\,+\,
x_3\,z_1z_2$$ generate the algebra $\textrm{Pol}(\fp_{\star})^K$. We can show that
the invariants $\alpha_1,\,\alpha_2,\,\xi$ and $\varphi$ are {\it algebraically independent}.
I will omit the detail.

\vskip 0.21cm
Now we will compute the $GL_{2,1}$-invariant differential
operators $D_1,\,D_2,\,\Psi,\,\Delta$ on ${\mathcal P}_{2,1}$
corresponding to the $K$-invariants $\alpha_1,\,\alpha_2,\,\xi,\,\varphi$
respectively under a canonical linear bijection
$$\Theta_{2,1}\,:\,
\textrm{Pol}(\fp_{\star})^K\lrt {\mathbb D}({\mathcal P}_{2,1}).$$
For real variables $t=(t_1,t_2,t_3)$ and $s=(s_1,s_2)$, we have
\begin{eqnarray*}
 &\ &
\ \text{exp}\,(t_1e_1+t_2e_2+t_3e_3+s_1f_1+s_2f_2)\\
&\ & =\left(\begin{pmatrix} a_1(t,s) & a_3(t,s)\\ a_3(t,s) &
a_2(t,s)\end{pmatrix},
\,\left(\,b_1(t,s),\,b_2(t,s)\,\right)\,\right),
\end{eqnarray*}
where
\begin{eqnarray*}
 a_1(t,s)&=&1+t_1+{1\over {2!}}(t_1^2+t_3^2)+{1\over
{3!}}(t_1^3+2t_1t_3^2 +t_2t_3^2)+\cdots,\\
a_2(t,s)&=&1+t_2+{1\over {2!}}(t_2^2+t_3^2)+{1\over
{3!}}(t_1t_3^2+2t_2t_3^2+ t_2^3)+\cdots,\\
a_3(t,s)&=&t_3+{1\over {2!}}(t_1+t_2)t_3+{1\over
{3!}}(t_1t_2+ t_1^2+t_2^2+t_3^2)t_3+\cdots,\\
b_1(t,s)&=& s_1-{1\over {2!}}(s_1t_1+s_2t_3)+{1\over
{3!}}\left\{\,
s_1(t_1^2+t_3^2)+s_2(t_1t_3+t_2t_3)\right\}-\cdots,\\
b_2(t,s)&=&s_2-{1\over {2!}}(s_1t_3+s_2t_2)+{1\over {3!}}
\left\{\,s_1(t_1+t_2)t_3+s_2(t_2^2+t_3^2)\right\}-\cdots.
\end{eqnarray*}
For brevity, we write $a_i,\,b_k$ for $a_i(t,s),\,b_k(t,s)\
(\,i=1,2,3,\ k=1,2\,)$ respectively. We now fix an element
$(g,c)\in GL_{2,1}$ and write $$g=\begin{pmatrix} g_1 & g_{12} \\
g_{21} & g_2 \end{pmatrix}\ \ \ \text{and}\ \ \ c=(c_1,c_2).$$ We
put
$$\big(Y(t,s),\,V(t,s)\big)=\left(\,(g,c)\cdot \text{exp}\left(
\sum_{i=1}^3t_ie_i+\sum_{k=1}^2s_kf_k\right)\,\right)\cdot
(I_2,0)$$ with $$Y(t,s)=\begin{pmatrix} y_1(t,s) & y_3(t,s)\\
y_3(t,s) & y_2(t,s)\end{pmatrix}\ \ \ \text{and}\ \ \
V(t,s)=(v_1(t,s),\,v_2(t,s)).$$ By an easy computation, we obtain

\begin{eqnarray*}
 y_1&=&(g_1a_1+g_{12}a_3)^2+(g_1a_3+g_{12}a_2)^2,\\
y_2&=&(g_{21}a_1+g_2a_3)^2+(g_{21}a_3+g_2a_2)^2,\\
y_3&=&(g_1a_1+g_{12}a_3)(g_{21}a_1+g_2a_3)+(g_1a_3+g_{12}a_2)
(g_{21}a_3+g_2a_2),\\
v_1&=&(c_1+b_1a_1+b_2a_3)g_1+(c_2+b_1a_3+b_2a_2)g_{12},\\
v_2&=&(c_1+b_1a_1+b_2a_3)g_{21}+(c_2+b_1a_3+b_2a_2)g_2.
\end{eqnarray*}

Using the chain rule, we can easily compute the
$GL_{2,1}$-invariant differential operators $D_1=\Theta_{2,1}(\alpha_1),\
D_2=\Theta_{2,1}(\alpha_2),\, \Psi=\Theta_{2,1}(\xi)$ and
$\Delta=\,\Theta_{2,1}(\varphi)$.
They are given by
$$D_1\,=\,2\, \textrm{tr} \left(Y\Yd\right) =2\left(\,y_1{{\partial\ }\over
{\partial y_1}} +y_2{{\partial\ }\over {\partial
y_2}}+y_3{{\partial}\over {\partial y_3}}\right),$$

\begin{eqnarray*}
D_2&=&  \textrm{tr} \left(\left(\,2Y\Yd\right)^2\right)\\
&=& 3\,D_1\,+\,8\left(\,y_3^2\,{{\partial^2\ \ \ }\over {\partial
y_1\partial y_2}}\,+\,y_1y_3\,{{\partial^2\quad}\over {\partial y_1
\partial y_3}}\,+\,
y_2y_3\,{{\partial^2\quad}\over {\partial y_2\partial y_3}}\,\right)\\
&\ &\ +\, 4\,\left\{\,y_1^2\,{{\partial^2}\over {\partial
y_1^2}}\,+\, y_2^2\,{{\partial^2}\over {\partial
y_2^2}}\,+\,{\frac 12}\big(y_1y_2+y_3^2\big) \,{{\partial^2}\over
{\partial y_3^2}}\,\right\},\\
\Psi&=&  \textrm{tr} \left(\,Y\,{}^{{}^{{}^{{}^\text{\scriptsize $t$}}}}\!\!\!\left(\Vd\right)\left(\Vd\right)
\right)\\
&=& y_1\,{{\partial^2}\over {\partial
v_1^2}}\,+\,2y_3\,{{\partial^2}\over {\partial v_1\partial
v_2}}\,+\,y_2\,{{\partial^2}\over {\partial v_2^2}}
\end{eqnarray*}
and
\begin{eqnarray*}
\Delta&=& {{\partial}\over{\partial V}}\left(\,
2Y{{\partial}\over{\partial
Y}}\,\right)Y\,{}^{{}^{{}^{{}^\text{\scriptsize $t$}}}}\!\!\!\left({{\partial}\over
{\partial V}}\right)\\
&=& 2\,\left(\,y_1^2\,{{\partial^3}\over{\partial y_1\partial
v_1^2}} \,+\,2\,y_1y_3\,{{\partial^3}\over {\partial y_1\partial
v_1\partial v_2}}\,
+\,y_3^2\,{{\partial^3}\over{\partial y_1\partial v_2^2}}\,\right)\\
&\ &\ +\,2\,\left(\,y_3^2\,{{\partial^3}\over{\partial y_2\partial
v_1^2}}\,+\, 2\,y_2\,y_3\,{{\partial^3}\over{\partial y_2\partial
v_1\partial v_2}}\,+\,
y_2^2\,{{\partial^3}\over{\partial y_2\partial v_2^2}}\,\right)\\
&\ &\ +\,2\, \left\{\,y_1 y_3\,{{\partial^3}\over{\partial
y_3\partial v_1^2}}\,+\, \big(\,y_1
y_2\,+\,y_3^2\,\big)\,{{\partial^3}\over{\partial y_3\partial v_1
\partial v_2}}\,+\,y_2 y_3\,{{\partial^3}\over{\partial y_3\partial v_2^2}}\,
\right\}\\
&\ &\ +\,3\,\left(\,y_1\,{{\partial^2}\over{\partial v_1^2}}\,+\,
2y_3\,{{\partial^2}\over{\partial v_1\partial v_2}}\,+\,
y_2\,{{\partial^2}\over{\partial v_2^2}}\,\right).
\end{eqnarray*}
Clearly $D_1$ commutes with $D_2$ but $\Psi$ does not commute with
$D_1$ and $D_2.$ Indeed, we have the following noncommutation
relations
\begin{eqnarray*}
[D_1,\,\Psi]&=&D_1\Psi-
\Psi D_1\\
&=& 2\,\Psi
\end{eqnarray*}
and
\begin{eqnarray*}
[D_2,\,\Psi]&=& D_2\Psi
-\Psi D_2\\
&=& 2\,(\,2\,D_1\,-1\,)\,\Psi \\
&\ &\ -\,8\,\text{det}\,(Y)\cdot\text{det}\,
\left({{\partial}\over {\partial Y}}+\,{}^{{}^{{}^{{}^\text{\scriptsize $t$}}}}\!\!\!\left( {{\partial}\over
{\partial V}}\right){{\partial}\over {\partial V}}
\right)\\
&\ &\
+\,8\,\text{det}\,(Y)\cdot\text{det}\,\left({{\partial}\over{\partial
Y}} \right)\,-\,4\,\big(y_1y_2+y_3^2\big)\, {{\partial^3}\over
{\partial y_3\partial v_1\partial v_2}}.
\end{eqnarray*}
Hence the algebra ${\mathbb D}({\mathcal P}_{2,1})$ is {\it not}
commutative.


\vskip 0.5cm
\noindent
{\bf 5.2.} {\bf The Case $n=2$ and $m=2$}
\vskip 0.3cm
In this case,
$$GL_{2,2}=GL(2,\BR)\ltimes
\BR^{(2,2)},\ \ K=O(2)\ \ \text{and} \ \ GL_{2,2}/K={\mathcal
P}_2\times \BR^{(2,2)}={\mathcal P}_{2,2}.$$ We see easily that
$$\fp_{\star}=\,\left\{\,(X,Z)\ \vert\ X=\,^tX\in \BR^{(2,2)},\ \ Z\in
\BR^{(2,2)} \,\right\}.$$
Let $O_2$ be the $2\times 2$ zero matrix.
We put
$$e_1=\left(\begin{pmatrix} 1 & 0\\ 0 & 0 \end{pmatrix},O_2\right),\ \ \ e_{2}=\left(\begin{pmatrix} 0 &
0\\ 0 & 1 \end{pmatrix},O_2\right),\ \ \ e_3=\left(\begin{pmatrix} 0 & {\frac 12} \\
{\frac 12}  & 0\end{pmatrix},O_2\right)$$
and $$f_1=\left( O_2,\begin{pmatrix} 1 & 0\\ 0 & 0 \end{pmatrix} \right),
\ \ f_2=\left(O_2,\begin{pmatrix} 0 & 1\\ 0 & 0 \end{pmatrix}\right),
\ \ f_3=\left(O_2,\begin{pmatrix} 0 & 0\\ 1 & 0 \end{pmatrix}\right), \\
f_4=\left(O_2,\begin{pmatrix} 0 & 0\\ 0 & 1 \end{pmatrix}\right).$$ Then
$\{\,e_1,\,e_2,\,e_3,\,f_1,\,f_2,\,f_3,\,f_4\,\}$ forms a basis for
$\fp_{\star}.$
We write for variables $(X,Z)\in \fp_{\star}$ by $$X=\begin{pmatrix} x_1 & {\frac 12}x_3 \\
{\frac 12}x_3 & x_2 \end{pmatrix}\ \ \ \text{and} \ \ \
Z=\begin{pmatrix} z_{11} & z_{12} \\ z_{21} & z_{22} \end{pmatrix}.$$
According to Theorem 3.3, the algebra $\textrm{Pol}(\fp_{\star})^K$ is generated by
the following polynomials

\begin{eqnarray*}
\alpha_1(X,Z)&=& \textrm{tr}   (X)=x_1\,+\,x_2,\\
\alpha_2(X,Z)&=& \textrm{tr} (X^2)=
x_1^2\,+\,x_2^2\,+\,{\frac 12}\,x_3^2,\\
\beta_{11}^{(0)}(X,Z)&=&(Z\,^t\!Z)_{11}=\,z_{11}^2\,+\,z_{12}^2,\\
\beta_{12}^{(0)}(X,Z)&=&(Z\,^t\!Z)_{12}=\,z_{11}z_{21}+z_{12}z_{22},\\
\beta_{22}^{(0)}(X,Z)&=&(Z\,^t\!Z)_{22}=\,z_{21}^2\,+\,z_{22}^2,\\
\beta_{11}^{(1)}(X,Z)&=&(ZX\,^t\!Z)_{11}=\,x_1z_{11}^2\,+\,x_2z_{12}^2+ x_3z_{11}z_{12},\\
\beta_{12}^{(1)}(X,Z)&=&(ZX\,^t\!Z)_{12}=\,x_1z_{11}z_{21}+x_2z_{12}z_{22}+
{\frac 12}x_3 (z_{11}z_{22}+z_{12}z_{21}),\\
\beta_{22}^{(1)}(X,Z)&=&(ZX\,^t\!Z)_{22}=\,x_1z_{21}^2\,+\,x_2z_{22}^2 +x_3 z_{21}z_{22}.\\
\end{eqnarray*}
We set
\begin{equation*}
\Delta_{ab}:=\det
\begin{pmatrix}
\beta_{11}^{(a)} & \beta_{12}^{(b)} \\
\beta_{12}^{(a)} & \beta_{22}^{(b)}
\end{pmatrix}
\qquad \textrm{for}\ a,b=0,1.
\end{equation*}

\noindent
Minoru Itoh showed that the following equation
\begin{equation}
\alpha_1\,\Delta_{00}-\Delta_{01}-\Delta_{10}=0
\end{equation}
\noindent holds.

\vskip 0.5cm We take a coordinate $(Y,V)$ in ${\mathcal P}_{2,2}$ that is
\begin{equation*}
Y=\begin{pmatrix} y_1 & y_3 \\
y_3 & y_2
\end{pmatrix}
\quad \textrm{and}\quad
V=\begin{pmatrix}
v_{11} & v_{12}  \\ v_{21} & v_{22}
\end{pmatrix}.
\end{equation*}
We put
\begin{equation*}
{{\partial\ \ }\over {\partial Y}}=
\begin{pmatrix} {{\partial\ \ }\over {\partial y_1}} & {\frac 12}{{\partial\ \ }\over {\partial y_3}}\\
{\frac 12}{{\partial\ \ }\over {\partial y_3}} & {{\partial\ \ }\over {\partial y_2}}
\end{pmatrix}
\quad \textrm{and}\quad
{{\partial\ \ } \over {\partial V}}=
\begin{pmatrix}
{{\partial\quad}\over {\partial v_{11}}} & {{\partial\quad}\over {\partial v_{12}}}\\
{{\partial\quad}\over {\partial v_{21}}} & {{\partial\quad}\over {\partial v_{22}}}
\end{pmatrix}.
\end{equation*}
We consider the following differential operators
\begin{equation*}
D_i:=\,\textrm{tr} \left(\left( 2\,Y {{\partial\ \ }\over {\partial Y}} \right)^i\right),\quad i=1,2
\end{equation*}
and
\begin{equation*}
\Omega_{pq}^{(k)}=\,\left\{\,{{\partial}\over{\partial V}}\left(
2Y{{\partial\ \ }\over{\partial Y}}\right)^kY\,{}^{{}^{{}^{{}^\text{\scriptsize $t$}}}}\!\!\!\left({{\partial\ \ }\over
{\partial V}}\right)\,\right\}_{pq},\quad k=0,1, \ 1\leq p\leq q\leq 2.
\end{equation*}
We note that $D_1,\,D_2,\,\Omega_{11}^{(0)},\cdots ,\,\Omega_{22}^{(1)}$ are $GL_{2,2}$-invariant.
For brevity we put
$$\partial_{ij}= {{\partial\quad}\over {\partial v_{ij} }},\quad i,j=1,2.$$
It is easily seen that
\begin{eqnarray*}
 D_1 &=& \textrm{tr} \left( 2\,Y {{\partial\ \ }\over {\partial Y}}\right) = 2\,\sum_{i=1}^{3}
\, y_i {{\partial\ \ }\over {\partial y_i}},\\
D_2&=&  3\,D_1\,+\,8\left(\,y_3^2\,{{\partial^2\ \ \ }\over {\partial
y_1\partial y_2}}\,+\,y_1y_3\,{{\partial^2\quad}\over {\partial y_1
\partial y_3}}\,+\,
y_2y_3\,{{\partial^2\quad}\over {\partial y_2\partial y_3}}\,\right)\\
&\ &\ +\, 4\,\left\{\,y_1^2\,{{\partial^2}\over {\partial
y_1^2}}\,+\, y_2^2\,{{\partial^2}\over {\partial
y_2^2}}\,+\,{\frac 12}\big(y_1y_2+y_3^2\big) \,{{\partial^2}\over
{\partial y_3^2}}\,\right\},\\
\Omega_{11}^{(0)} &=& y_1\, \partial_{11}^2 + y_2 \,\partial_{12}^2
+ 2\,y_3\, \partial_{11}\partial_{12},\\
\Omega_{12}^{(0)} &=& y_1\, \partial_{11} \partial_{21}+ y_2 \,\partial_{12} \partial_{22}
+ y_3\,\left(  \partial_{11}\partial_{22} + \partial_{12}\partial_{21}\right),\\
\Omega_{22}^{(0)} &=& y_1\, \partial_{21}^2 + y_2 \,\partial_{22}^2
+ 2\,y_3\, \partial_{21}\partial_{22}.
\end{eqnarray*}
Then by a direct computation, we have the following relations
\begin{equation}
[D_1,D_2]=0,
\end{equation}
\begin{equation}
[\Omega_{kl}^{(0)},\Omega_{pq}^{(0)}]=0,\quad 1\leq k\leq l\leq 2,\ 1\leq p\leq q\leq 2,
\end{equation}
\begin{equation}
[D_1,\Omega_{11}^{(0)}]=2\,\Omega_{11}^{(0)},\quad [D_1,\Omega_{12}^{(0)}]=2\,\Omega_{12}^{(0)},\quad
[D_1,\Omega_{22}^{(0)}]=2\,\Omega_{22}^{(0)}.
\end{equation}
Therefore $\BD ({\mathcal P}_{2,2})$ is not commutative.

\end{section}

\vskip 0.87cm

\begin{section}{{\bf The Case $n=3$}}
\setcounter{equation}{0}

\vskip 0.5cm
\noindent
{\bf 6.1.} {\bf The Case $n=3$ and $m=1$}
\vskip 0.3cm
In this case,
$$GL_{3,1}=GL(3,\BR)\ltimes
\BR^{(1,3)},\ \ K=O(3)\ \ \text{and} \ \ GL_{3,1}/K={\mathcal
P}_3\times \BR^{(1,3)}={\mathcal P}_{3,1}.$$ We see easily that
$$\fp_{\star}=\,\left\{\,(X,Z)\ \vert\ X=\,^tX\in \BR^{(3,3)},\ \ Z\in
\BR^{(1,3)} \,\right\}.$$

\noindent We put
\begin{eqnarray*}
E_1&=&\begin{pmatrix} 1 & 0 & 0\\
0 & 0 & 0 \\ 0 & 0 & 0 \end{pmatrix},\quad
E_2=\begin{pmatrix} 0 & 0 & 0\\
0 & 1 & 0 \\ 0 & 0 & 0 \end{pmatrix},\quad
E_3=\begin{pmatrix} 0 & 0 & 0\\
0 & 0 & 0 \\ 0 & 0 & 1 \end{pmatrix},\\
E_4&=&\begin{pmatrix} 0 & {\frac 12} & 0\\
{\frac 12} & 0 & 0 \\ 0 & 0 & 0 \end{pmatrix},\quad
E_5=\begin{pmatrix} 0 & 0 & {\frac 12}\\
0 & 0 & 0 \\ {\frac 12} & 0 & 0 \end{pmatrix},\quad
E_6=\begin{pmatrix} 0 & 0 & 0\\
0 & 0 & {\frac 12} \\ 0 & {\frac 12} & 0 \end{pmatrix}.
\end{eqnarray*}

Let $O_3$ be the $3\times 3$ zero matrix and let $O_{1,3}=(0,0,0)\in\BR^{(1,3)}$.
We put
\begin{eqnarray*}
e_i&=& (E_i,O_{1,3}),\quad 1\leq i\leq 6,\\
f_1&=& (O_3,(1,0,0)),\quad f_2=(O_3,(0,1,0)),\quad f_3=(O_3,(0,0,1)).
\end{eqnarray*}
Then
$\{\,e_i,\,f_j\,|\,1\leq i\leq 6,\ 1\leq j\leq 3\,\}$ forms a basis for
$\fp_{\star}.$
Using this basis,
We write for variables $(X,Z)\in \fp_{\star}$ by
$$X=\begin{pmatrix} x_1 & {\frac 12}x_4 & {\frac 12}x_5 \\
{\frac 12}x_4 & x_2 & {\frac 12}x_6\\
{\frac 12}x_5 & {\frac 12}x_6 & x_3 \end{pmatrix}\ \ \ \text{and} \ \ \
Z=(z_1,z_2,z_3).$$

According to Theorem 3.3, the algebra $\textrm{Pol}(\fp_{\star})^K$ is generated by
the following polynomials

\begin{eqnarray*}
\alpha_1(X,Z) &=& x_1+x_2+x_3,\\
\alpha_2(X,Z) &=& x_1^2+x_2^2+x_3^2+ {\frac 12}\big(x_4^2+x_5^2+x_6^2\big),\\
\alpha_3(X,Z) &=& x_1^3+x_2^3+x_3^3+ {\frac 34}\,\left\{(x_1+x_2)x_4^2 + (x_1+x_3)x_5^2
+ (x_2+x_3)x_6^2\right\}\\
& &\,  +\,{\frac 34}\,x_4x_5x_6,\\
\beta_0(X,Z)&=& z_1^2+z_2^2+z_3^2,\\
\beta_1(X,Z)&=& x_1z_1^2+x_2z_2^2+x_3z_3^2 +x_4z_1z_2+x_5z_1z_3+x_6z_2z_3,\\
\beta_2(X,Z)&=& x_1^2z_1^2+x_2^2z_2^2+\,{\frac 14}\,
\left\{ \big(x_4^2+x_5^2\big)z_1^2\,+\,\big(x_4^2+x_6^2\big)z_2^2\,+\,
\big(x_5^2+x_6^2\big)z_3^2\right\}\\
& &\, +\,\left( x_1x_4+x_2x_4+ {\frac 12}x_5x_6\right)z_1z_2\,+\,
\left( x_1x_5+x_3x_5+ {\frac 12}x_4x_6\right)z_1z_3\\
& &\, +\,\left( x_2x_6+x_3x_6+ {\frac 12}x_4x_5\right)z_2z_3.
\end{eqnarray*}
M. Itoh claims that $\alpha_1,\alpha_2,\alpha_3,\beta_0,\beta_1,\beta_2$ are {\it algebraically independent}.

\vskip 0.5cm We take a coordinate $(Y,V)$ in ${\mathcal P}_{3,1}$ that is
\begin{equation*}
Y=\begin{pmatrix} y_1 & y_4 & y_5 \\
y_4 & y_2 & y_6 \\
y_5 & y_6 & y_3
\end{pmatrix}
\quad \textrm{and}\quad
V=(v_1,v_2,v_3).
\end{equation*}
We put
\begin{equation*}
{{\partial\ \ }\over {\partial Y}}=
\begin{pmatrix} {{\partial\ \ }\over {\partial y_1}} & {\frac 12}{{\partial\ \ }\over {\partial y_4}}
 & {\frac 12}{{\partial\ \ }\over {\partial y_5}} \\
{\frac 12}{{\partial\ \ }\over {\partial y_4}} & {{\partial\ \ }\over {\partial y_2}}
& {\frac 12}{{\partial\ \ }\over {\partial y_6}} \\
{\frac 12}{{\partial\ \ }\over {\partial y_5}} & {\frac 12}{{\partial\ \ }\over {\partial y_6}}
& {{\partial\ \ }\over {\partial y_3}}
\end{pmatrix}
\quad \textrm{and}\quad
{{\partial\ \ } \over {\partial V}}=
\left( {{\partial\quad}\over {\partial v_{1}}}, {{\partial\quad}\over {\partial v_{2}}},
{{\partial\quad}\over {\partial v_{3}}}\right).
\end{equation*}
We consider the following differential operators
\begin{equation*}
D_i:=\,\textrm{tr} \left(\left( 2\,Y {{\partial\ \ }\over {\partial Y}} \right)^i\right),\quad i=1,2,3
\end{equation*}
and
\begin{equation*}
\Omega_k=\,{{\partial}\over{\partial V}}\left(
2Y{{\partial\ \ }\over{\partial Y}}\right)^kY\,{}^{{}^{{}^{{}^\text{\scriptsize $t$}}}}\!\!\!\left({{\partial\ \ }\over
{\partial V}}\right),\quad k=0,1,2.
\end{equation*}
We note that $D_1,\,D_2,\,D_3,\,\Omega_0,\,\Omega_1$ and $\Omega_2$ are $GL_{2,2}$-invariant.
It is easily seen that
\begin{eqnarray*}
 D_1 &=& \textrm{tr} \left( 2\,Y {{\partial\ \ }\over {\partial Y}}\right) = 2\,\sum_{i=1}^{6}
\, y_i {{\partial\ \ }\over {\partial y_i}},\\
\Omega_0&=& y_1\ {{\partial^2\ \ }\over {\partial v_1^2}} + y_2\,{{\partial^2\ \ }\over {\partial v_2^2}}
+ y_3\,{{\partial^2\ \ }\over {\partial v_3^2}}\\
& & +\, 2\,y_4\, {{\partial^2\quad\ }\over {\partial v_1\partial v_2}}
+\, 2\,y_5 \, {{\partial^2\quad\ }\over {\partial v_1\partial v_3}}
+\, 2\,y_6 \, {{\partial^2\quad\ }\over {\partial v_2\partial v_3}}.
\end{eqnarray*}
\noindent
Then we have the following relations
\begin{equation}
[D_i,D_j]=0 \quad \textrm{for all}\ i,j=1,2,3
\end{equation}
and
\begin{equation}
[D_1,\Omega_0]=2\,\Omega_0.
\end{equation}
Therefore $\BD ({\mathcal P}_{3,1})$ is not commutative.


\vskip 0.5cm
\noindent
{\bf 6.2.} {\bf The Case $n=3$ and $m=2$}
\vskip 0.3cm
In this case,
$$GL_{3,2}=GL(3,\BR)\ltimes
\BR^{(2,3)},\ \ K=O(3)\ \ \text{and} \ \ GL_{3,2}/K={\mathcal
P}_3\times \BR^{(2,3)}={\mathcal P}_{3,2}.$$ We see easily that
$$\fp_{\star}=\,\left\{\,(X,Z)\ \vert\ X=\,^tX\in \BR^{(3,3)},\ \ Z\in
\BR^{(2,3)} \,\right\}.$$

\noindent We put
\begin{eqnarray*}
E_1&=&\begin{pmatrix} 1 & 0 & 0\\
0 & 0 & 0 \\ 0 & 0 & 0 \end{pmatrix},\quad
E_2=\begin{pmatrix} 0 & 0 & 0\\
0 & 1 & 0 \\ 0 & 0 & 0 \end{pmatrix},\quad
E_3=\begin{pmatrix} 0 & 0 & 0\\
0 & 0 & 0 \\ 0 & 0 & 1 \end{pmatrix},\\
E_4&=&\begin{pmatrix} 0 & {\frac 12} & 0\\
{\frac 12} & 0 & 0 \\ 0 & 0 & 0 \end{pmatrix},\quad
E_5=\begin{pmatrix} 0 & 0 & {\frac 12}\\
0 & 0 & 0 \\ {\frac 12} & 0 & 0 \end{pmatrix},\quad
E_6=\begin{pmatrix} 0 & 0 & 0\\
0 & 0 & {\frac 12} \\ 0 & {\frac 12} & 0 \end{pmatrix}.
\end{eqnarray*}
and
\begin{eqnarray*}
F_1&=&\begin{pmatrix} 1 & 0 & 0\\
0 & 0 & 0  \end{pmatrix},\quad
F_2=\begin{pmatrix} 0 & 1 & 0\\
0 & 0 & 0  \end{pmatrix},\quad
F_3=\begin{pmatrix} 0 & 0 & 1\\
0 & 0 & 0  \end{pmatrix},\\
F_4&=&\begin{pmatrix} 0 & 0 & 0\\
1 & 0 & 0  \end{pmatrix},\quad
F_5=\begin{pmatrix} 0 & 0 & 0\\
0 & 1 & 0 \end{pmatrix},\quad
F_6=\begin{pmatrix} 0 & 0 & 0\\
0 & 0 & 1 \end{pmatrix}.
\end{eqnarray*}

Let $O_3$ be the $3\times 3$ zero matrix and let
$$O_{2,3}=\begin{pmatrix} 0 & 0 & 0 \\ 0 & 0 & 0 \end{pmatrix} \in\BR^{(2,3)}$$.
We put
$$e_i=(E_i,O_{2,3}),\quad f_j=(O_3,F_j)\quad 1\leq i,j\leq 6.$$
Then
$\{\,e_i,\,f_j\,|\,1\leq i,j\leq 6\,\}$ forms a basis for
$\fp_{\star}.$
Using this basis,
We write for variables $(X,Z)\in \fp_{\star}$ by
$$X=\begin{pmatrix} x_1 & {\frac 12}x_4 & {\frac 12}x_5 \\
{\frac 12}x_4 & x_2 & {\frac 12}x_6\\
{\frac 12}x_5 & {\frac 12}x_6 & x_3 \end{pmatrix}\ \ \ \text{and} \ \ \
Z=\begin{pmatrix} z_{11} & z_{12} & z_{13}\\
z_{21} & z_{22} & z_{23} \end{pmatrix}.$$

According to Theorem 3.3, the algebra $\textrm{Pol}(\fp_{\star})^K$ is generated by
the following polynomials

\begin{eqnarray*}
\alpha_1(X,Z) &=& x_1+x_2+x_3,\\
\alpha_2(X,Z) &=& x_1^2+x_2^2+x_3^2+ {\frac 12}\big(x_4^2+x_5^2+x_6^2\big),\\
\alpha_3(X,Z) &=& x_1^3+x_2^3+x_3^3+ {\frac 34}\,\left\{(x_1+x_2)x_4^2 + (x_1+x_3)x_5^2
+ (x_2+x_3)x_6^2 \right\}\\
& &\,  +\,{\frac 34}\,x_4x_5x_6,\\
\beta_{11}^{(0)}(X,Z)&=& z_{11}^2+z_{12}^2+z_{13}^2,\\
\beta_{12}^{(0)}(X,Z)&=& z_{11}z_{21}+z_{12}z_{22}+z_{13}z_{23},\\
\beta_{22}^{(0)}(X,Z)&=& z_{21}^2+z_{22}^2+z_{23}^2,\\
\beta_{11}^{(1)}(X,Z)&=& x_1z_{11}^2+x_2z_{12}^2+x_3z_{13}^2 +x_4z_{11}z_{12}+x_5z_{11}z_{13}+x_6z_{12}z_{13},\\
\beta_{12}^{(1)}(X,Z)&=& x_1z_{11}z_{21}+x_2z_{12}z_{22}+x_3z_{13}z_{23}\,+\,
{\frac 12}\,x_4 (z_{11}z_{22}+z_{12}z_{21})\\
& & \, +\,{\frac 12}\,x_5 (z_{11}z_{23}+z_{13}z_{21})\,+\,{\frac 12}\,x_6 (z_{12}z_{23}+z_{13}z_{22}),\\
\beta_{22}^{(1)}(X,Z)&=&  x_1z_{21}^2+x_2z_{22}^2+x_3z_{23}^2 +x_4z_{21}z_{22}+x_5z_{21}z_{23}+x_6z_{22}z_{23},  \\
\beta_{11}^{(2)}(X,Z)&=& x_1^2z_{11}^2+x_2^2z_{12}^2+ x_3^2z_{13}^2 \\
& &\,+\, {\frac 14}\,\left\{\, x_4^2 \big( z_{11}^2+z_{12}^2\big)\,+\,x_5^2 \big( z_{11}^2+z_{13}^2\big)\,+\,
x_6^2 \big( z_{12}^2+z_{13}^2\big)\,\right\}\\
& &\, +\,(x_1+x_2)\,x_4z_{11}z_{12}\,+\,(x_1+x_3)\,x_5z_{11}z_{13}\,+\,(x_2+x_3)\,x_6z_{12}z_{13}\\
& &\,+\, {\frac 12} \,\big( x_4x_5z_{12}z_{13}\,+\,x_4x_6z_{11}z_{13}\,+\,x_5x_6z_{11}z_{12}\big),   \\
\beta_{12}^{(2)}(X,Z)&=&  x_1^2z_{11}z_{21}+x_2^2z_{12}z_{22}+ x_3^2z_{13}z_{23} \\
& &\,+\, {\frac 14}\,\left\{ \big( x_4^2+x_5^2\big)\,z_{11}z_{21}+
\big( x_4^2+x_6^2\big)\,z_{12}z_{22}+
\big( x_5^2+x_6^2\big)\,z_{13}z_{23} \right\}
\\
& &\,+\,{\frac 12}\,\left( x_1x_4+x_2x_4+{\frac 12}x_5x_6\right)\,(z_{11}z_{22}+z_{12}z_{21})\\
& &\,+\,{\frac 12}\,\left( x_1x_5+x_3x_5+{\frac 12}x_4x_6\right)\,(z_{11}z_{23}+z_{13}z_{21})\\
& &\,+\,{\frac 12}\,\left( x_2x_6+x_3x_6+{\frac 12}x_4x_5\right)\,(z_{12}z_{23}+z_{13}z_{22}),\\
\beta_{22}^{(2)}(X,Z)&=& x_1^2z_{21}^2+x_2^2z_{22}^2+ x_3^2z_{23}^2 \\
& &\,+\, {\frac 14}\,\left\{\, x_4^2 \big( z_{21}^2+z_{22}^2\big)\,+\,x_5^2 \big( z_{21}^2+z_{23}^2\big)\,+\,
x_6^2 \big( z_{22}^2+z_{23}^2\big)\,\right\}\\
& &\, +\,(x_1+x_2)\,x_4z_{21}z_{22}\,+\,(x_1+x_3)\,x_5z_{21}z_{23}\,+\,(x_2+x_3)\,x_6z_{22}z_{23}\\
& &\,+\, {\frac 12} \,\big( x_4x_5z_{22}z_{23}\,+\,x_4x_6z_{21}z_{23}\,+\,x_5x_6z_{21}z_{22}\big).
\end{eqnarray*}

We set
\begin{equation*}
\Delta_{ab}:=\det
\begin{pmatrix}
\beta_{11}^{(a)} & \beta_{12}^{(b)} \\
\beta_{12}^{(a)} & \beta_{22}^{(b)}
\end{pmatrix}
\qquad \textrm{for}\ a,b=0,1,2.
\end{equation*}

Minoru Itoh claims that
\begin{equation}
\big( \alpha_1^2-\alpha_2\big)\,\Delta_{00} - 2\,\alpha_1 \,\big( \Delta_{01}+\Delta_{10}\big)
+\,2\,\big( \Delta_{02}+\Delta_{11}+\Delta_{20}\big)=0.
\end{equation}

\vskip 0.5cm We take a coordinate $(Y,V)$ in ${\mathcal P}_{3,2}$ that is
\begin{equation*}
Y=\begin{pmatrix} y_1 & y_4 & y_5 \\
y_4 & y_2 & y_6 \\
y_5 & y_6 & y_3
\end{pmatrix}
\quad \textrm{and}\quad
V=\begin{pmatrix}
v_{11} & v_{12} & v_{13}  \\ v_{21} & v_{22} & v_{23}
\end{pmatrix}.
\end{equation*}
We put
\begin{equation*}
{{\partial\ \ }\over {\partial Y}}=
\begin{pmatrix} {{\partial\ \ }\over {\partial y_1}} & {\frac 12}{{\partial\ \ }\over {\partial y_4}}
 & {\frac 12}{{\partial\ \ }\over {\partial y_5}} \\
{\frac 12}{{\partial\ \ }\over {\partial y_4}} & {{\partial\ \ }\over {\partial y_2}}
& {\frac 12}{{\partial\ \ }\over {\partial y_6}} \\
{\frac 12}{{\partial\ \ }\over {\partial y_5}} & {\frac 12}{{\partial\ \ }\over {\partial y_6}}
& {{\partial\ \ }\over {\partial y_3}}
\end{pmatrix}
\quad \textrm{and}\quad
\begin{pmatrix}
{{\partial\quad}\over {\partial v_{11}}} & {{\partial\quad}\over {\partial v_{12}}}
& {{\partial\quad}\over {\partial v_{13}}}  \\
{{\partial\quad}\over {\partial v_{21}}} & {{\partial\quad}\over {\partial v_{22}}}
& {{\partial\quad}\over {\partial v_{23}}}
\end{pmatrix}.
\end{equation*}
\vskip 0.3cm\noindent
We consider the following differential operators
\begin{equation*}
D_i:=\,\textrm{tr} \left(\left( 2\,Y {{\partial\ \ }\over {\partial Y}} \right)^i\right),\quad i=1,2,3
\end{equation*}
and
\begin{equation*}
\Omega_{pq}^{(k)}=\,\left\{\,{{\partial}\over{\partial V}}\left(
2Y{{\partial\ \ }\over{\partial Y}}\right)^kY\,{}^{{}^{{}^{{}^\text{\scriptsize $t$}}}}\!\!\!\left({{\partial\ \ }\over
{\partial V}}\right)\,\right\}_{pq},\quad k=0,1,2, \ 1\leq p\leq q\leq 2.
\end{equation*}
We note that $D_1,\,D_2,\,D_3,\,\Omega_{11}^{(0)},\cdots ,\,\Omega_{22}^{(2)}$ are $GL_{3,2}$-invariant.
For brevity we put
$$\partial_{ij}= {{\partial\quad}\over {\partial v_{ij} }},\quad i=1,2,\ j=1,2,3.$$
It is easily seen that
\begin{eqnarray*}
 D_1 &=& \textrm{tr} \left( 2\,Y {{\partial\ \ }\over {\partial Y}}\right) = 2\,\sum_{i=1}^{6}
\, y_i {{\partial\ \ }\over {\partial y_i}},\\
\Omega_{11}^{(0)}&=& y_1\,\partial_{11}^2 +\,y_2\,\partial_{12}^2 +\,y_3\,\partial_{13}^2
 +2\,y_4\,\partial_{11}\partial_{12} +2\,y_5\,\partial_{11}\partial_{13} +2\,y_6\,\partial_{12}\partial_{13}, \\
\Omega_{12}^{(0)}&=& y_1\,\partial_{11}\partial_{21} + y_2\,\partial_{12}\partial_{22}
+ y_3\,\partial_{13}\partial_{23} + y_4\,\big( \partial_{11}\partial_{22}+\partial_{12}\partial_{21}\big) \\
& & + y_5\,\big( \partial_{11}\partial_{23}+\partial_{13}\partial_{21}\big)
+ y_6\,\big( \partial_{12}\partial_{23}+\partial_{13}\partial_{22}\big),\\
\Omega_{22}^{(0)}&=& y_1\,\partial_{21}^2 +\,y_2\,\partial_{22}^2 +\,y_3\,\partial_{23}^2
 +2\,y_4\,\partial_{21}\partial_{22} +2\,y_5\,\partial_{21}\partial_{23} +2\,y_6\,\partial_{22}\partial_{23}.
\end{eqnarray*}
\noindent
Then we have the following relations
\begin{equation}
[D_i,D_j]=0 \quad \textrm{for all}\ i,j=1,2,3,
\end{equation}
\begin{equation}
[\Omega_{kl}^{(0)},\Omega_{pq}^{(0)}]=0,\quad 1\leq k\leq l\leq 2,\ 1\leq p\leq q\leq 2
\end{equation}
and
\begin{equation}
[D_1,\Omega_{11}^{(0)}]=2\,\Omega_{11}^{(0)},\quad [D_1,\Omega_{12}^{(0)}]=2\,\Omega_{12}^{(0)},\quad
[D_1,\Omega_{22}^{(0)}]=2\,\Omega_{22}^{(0)}.
\end{equation}
Therefore $\BD ({\mathcal P}_{3,2})$ is not commutative.

\end{section}

\vskip 0.87cm


\begin{section}{{\bf The Case $n=4$}}
\setcounter{equation}{0}

\vskip 0.5cm
\noindent
{\bf 6.1.} {\bf The Case $n=4$ and $m=1$}
\vskip 0.3cm
In this case,
$$GL_{4,1}=GL(4,\BR)\ltimes
\BR^{(1,4)},\ \ K=O(4)\ \ \text{and} \ \ GL_{4,1}/K={\mathcal
P}_4\times \BR^{(1,4)}={\mathcal P}_{4,1}.$$ We see easily that
$$\fp_{\star}=\,\left\{\,(X,Z)\ \vert\ X=\,^tX\in \BR^{(4,4)},\ \ Z\in
\BR^{(1,4)} \,\right\}.$$

\noindent We put
\begin{eqnarray*}
E_1&=&\begin{pmatrix} 1 & 0 & 0 & 0\\
0 & 0 & 0  & 0\\ 0 & 0 & 0  & 0 \\ 0 & 0 & 0  & 0 \end{pmatrix},\
E_2=\begin{pmatrix} 0 & 0 & 0 & 0\\
0 & 1 & 0  & 0\\ 0 & 0 & 0  & 0 \\ 0 & 0 & 0  & 0 \end{pmatrix},\
E_3=\begin{pmatrix} 0 & 0 & 0 & 0\\
0 & 0 & 0  & 0\\ 0 & 0 & 1  & 0 \\ 0 & 0 & 0  & 0 \end{pmatrix},\\
E_4&=&\begin{pmatrix} 0 & 0 & 0 & 0\\
0 & 0 & 0  & 0\\ 0 & 0 & 0  & 0 \\ 0 & 0 & 0  & 1 \end{pmatrix},
\ E_5=\begin{pmatrix} 0 & {\frac 12} & 0 & 0\\
{\frac 12} & 0 & 0  & 0\\ 0 & 0 & 0  & 0 \\ 0 & 0 & 0  & 0 \end{pmatrix},
\ \ E_6=\begin{pmatrix} 0 & 0 & {\frac 12} & 0\\
0 & 0 & 0  & 0\\ {\frac 12} & 0 & 0  & 0 \\ 0 & 0 & 0  & 0 \end{pmatrix},
\
\\
E_7&=&\begin{pmatrix} 0 & 0 & 0 & {\frac 12} \\
0 & 0 & 0  & 0\\ 0 & 0 & 0  & 0 \\ {\frac 12} & 0 & 0  & 0 \end{pmatrix},\
E_8=\begin{pmatrix} 0 & 0 & 0 & 0 \\
0 & 0 & {\frac 12}  & 0\\ 0 &  {\frac 12}& 0  & 0 \\ 0 & 0 & 0  & 0 \end{pmatrix},\\
E_9&=& \begin{pmatrix} 0 & 0 & 0 & 0 \\
0 & 0 & 0 & {\frac 12}  \\ 0 & 0 & 0  & 0 \\ 0 & {\frac 12} & 0  & 0 \end{pmatrix},\
E_{10}=\begin{pmatrix} 0 & 0 & 0 & 0 \\
0 & 0 & 0 & 0  \\ 0 & 0 & 0  & {\frac 12}  \\ 0 & 0 & {\frac 12} & 0   \end{pmatrix}.
\end{eqnarray*}

Let $O_4$ be the $4\times 4$ zero matrix and let $O_{1,4}=(0,0,0,0)\in\BR^{(1,4)}$.
We put
\begin{eqnarray*}
e_i&=& (E_i,O_{1,4}),\quad 1\leq i\leq 10,\\
f_1&=& (O_4,(1,0,0,0)),\ f_2=(O_4,(0,1,0,0)),\\
f_3&=&(O_4,(0,0,1,0)),\ f_4=(O_4,(0,0,0,1)).
\end{eqnarray*}
Then
$\{\,e_i,\,f_j\,|\,1\leq i\leq 10,\ 1\leq j\leq 4\,\}$ forms a basis for
$\fp_{\star}.$
Using this basis,
We write for variables $(X,Z)\in \fp_{\star}$ by
$$X=\begin{pmatrix} x_1 & {\frac 12}x_5 & {\frac 12}x_6  & {\frac 12}x_7\\
{\frac 12}x_5 & x_2 & {\frac 12}x_8  & {\frac 12}x_9\\
{\frac 12}x_6 & {\frac 12}x_8 &x_3 & {\frac 12}x_{10}\\
{\frac 12}x_7 & {\frac 12}x_9 & {\frac 12}x_{10} & x_4
\end{pmatrix}
 \ \ \text{and} \ \
Z=(z_1,z_2,z_3,z_4).$$
We put
\begin{eqnarray}
A&=& x_1^2+{\frac 14}\,x_5^2+{\frac 14}\,x_6+{\frac 14}\,x_7^2,\\
B&=& x_2^2+{\frac 14}\,x_5^2+{\frac 14}\,x_8+{\frac 14}\,x_9^2,\\
C&=& x_3^2+{\frac 14}\,x_6^2+{\frac 14}\,x_8+{\frac 14}\,x_{10}^2,\\
D&=& x_4^2+{\frac 14}\,x_7^2+{\frac 14}\,x_9+{\frac 14}\,x_{10}^2,\\
E&=& {\frac 12}\,(x_1+x_2)\,x_5+{\frac 14}\,(x_6\,x_8+x_7\,x_9),\\
F&=& {\frac 12}\,(x_1+x_3)\,x_6+{\frac 14}\,(x_3\,x_6+x_5\,x_8),\\
G&=& {\frac 12}\,(x_1+x_4)\,x_7+{\frac 14}\,(x_5\,x_9+x_6\,x_{10}),\\
H&=& {\frac 12}\,(x_2+x_3)\,x_8+{\frac 14}\,(x_5\,x_6+x_9\,x_{10}),\\
I&=& {\frac 12}\,(x_2+x_4)\,x_9+{\frac 14}\,(x_5\,x_7+x_8\,x_{10}),\\
J&=& {\frac 12}\,(x_3+x_4)\,x_{10}+{\frac 14}\,(x_6\,x_{10}+x_6\,x_7).
\end{eqnarray}

According to Theorem 3.3, the algebra $\textrm{Pol}(\fp_{\star})^K$ is generated by
the following polynomials

\begin{eqnarray*}
\alpha_1(X,Z)&=& x_1+x_2+x_3+x_4,\\
\alpha_2(X,Z)&=& x_1^2+x_2^2+x_3^2+x_4^2+{\frac 12}\,\big(x_5^2+x_6^2+x_7^2+x_8^2+x_9^2+x_{10}^2\big),\\
\alpha_3(X,Z)&=& x_1^3+x_2^3+x_3^3+x_4^3\\
& & + {\frac 34}\,x_1 \big( x_5^2+x_6^2+x_7^2\big)+ {\frac 34}\,x_2 \big( x_5^2+x_8^2+x_9^2\big)  \\
& & + {\frac 34}\,x_3 \big( x_6^2+x_8^2+x_{10}^2\big)+ {\frac 34}\,x_4 \big( x_7^2+x_9^2+x_{10}^2\big)  \\
& & + {\frac 34}\,\big( x_5\, x_6\, x_8 + x_5 \, x_7\, x_9 + x_6\, x_7\, x_{10} + x_8\, x_9 \, x_{10} \big),
\end{eqnarray*}

\newpage

\begin{eqnarray*}
\alpha_4(X,Z)&=& A^2+B^2+C^2+D^2+ 2\,\big( E^2+F^2+G^2+H^2+I^2+J^2\big),\\
\beta_0(X,Z)&=& z_1^2 + z_2^2 + z_3^2 + z_4^2,\\
\beta_1(X,Z)&=& x_1z_1^2+x_2z_2^2+x_3z_3^2+x_4z_4^2,\\
& & + x_5 z_1 z_2 + x_6 z_1 z_3 + x_7 z_1 z_4 + x_8 z_2 z_3 + x_9 z_2 z_4 + x_{10} z_3 z_4,\\
\beta_2(X,Z)&=& A z_1^2+ B z_2^2+ C z_3^2+ D z_4^2,\\
& & + 2 \big( E z_1 z_2 + F z_1 z_3 + G z_1 z_4 + H z_2 z_3 + I z_2 z_4 + J z_3 z_4 \big),\\
\beta_3(X,Z)&=& {\frac 12}\,\big( 2 Ax_1 + Ex_5 + Fx_6 + Gx_7 \big) z_1^2 +
{\frac 12}\,\big( 2 B x_2 + E x_5 + H x_8 + I x_9 \big) z_2^2 \\
& & + {\frac 12}\,\big( 2 C x_3 + F x_6 + H x_8 + J x_{10} \big) z_3^2 +
{\frac 12}\,\big( 2 D x_4 + G x_7 + I x_9 + J x_{10} \big) z_4^2 \\
& & + {\frac 12}\,\left\{ 2E(x_1+x_2)+(A+B)x_5+ H x_6 + I x_7 + F x_8 + G x_9 \right\} z_1z_2\\
& & + {\frac 12}\,\left\{ 2F(x_1+x_3)+ H x_5 + (A+C)x_6+  J x_7 + E x_8 + G x_{10} \right\} z_1z_3\\
& & + {\frac 12}\,\left\{ 2G(x_1+x_4)+ I x_5 + J x_6 + (A+D)x_7 +  E x_9 + F x_{10} \right\} z_1z_4\\
& & + {\frac 12}\,\left\{ 2H(x_2+x_3)+ F x_5 + E x_6 + (B+C)x_8 +  J x_9 + I x_{10} \right\} z_2z_3\\
& & + {\frac 12}\,\left\{ 2I(x_2+x_4)+ G x_5 + E x_7 + J x_8 + (B+D)x_9 +  H x_{10} \right\} z_2z_4\\
& & + {\frac 12}\,\left\{ 2J(x_3+x_4)+ G x_6 + F x_7 + I x_8 + H x_9 + (C+D) x_{10} \right\} z_3z_4.
\end{eqnarray*}
M. Itoh claims that the above eight invariants $\alpha_1,\cdots,\alpha_4,\,\beta_0,\cdots,\beta_3$
are {\it algebraically independent}.

\vskip 0.5cm We take a coordinate $(Y,V)$ in ${\mathcal P}_{4,1}$ that is
\begin{equation*}
Y=\begin{pmatrix} y_1 & y_5 & y_6 & y_7 \\ y_5 & y_2 & y_8 & y_9 \\
y_6 & y_8 & y_3 & y_{10} \\ y_7 & y_9 & y_{10} & y_4
\end{pmatrix}
\quad \textrm{and}\quad
V=(v_1,v_2,v_3,v_4).
\end{equation*}

We put
\begin{equation*}
{{\partial\ \ }\over {\partial Y}}=
\begin{pmatrix} {{\partial\ \ }\over {\partial y_1}} & {\frac 12}{{\partial\ \ }\over {\partial y_5}}
& {\frac 12}{{\partial\ \ }\over {\partial y_6}}  & {\frac 12}{{\partial\ \ }\over {\partial y_7}} \\
{\frac 12}{{\partial\ \ }\over {\partial y_5}} & {{\partial\ \ }\over {\partial y_2}}
& {\frac 12}{{\partial\ \ }\over {\partial y_8}} & {\frac 12}{{\partial\ \ }\over {\partial y_9}} \\
{\frac 12}{{\partial\ \ }\over {\partial y_6}} & {\frac 12}{{\partial\ \ }\over {\partial y_8}}
& {{\partial\ \ }\over {\partial y_3}} & {\frac 12}{{\partial\ \ }\over {\partial y_{10}}} \\
{\frac 12}{{\partial\ \ }\over {\partial y_7}} & {\frac 12}{{\partial\ \ }\over {\partial y_9}}
& {\frac 12} { {\partial\ \ }\over {\partial y_{10}} } & { {\partial\ \ }\over {\partial y_4} }
\end{pmatrix}
\quad \textrm{and}\quad
{{\partial\ \ } \over {\partial V}}=
\left( {{\partial\ \ }\over {\partial v_{1}}},  {{\partial\ \ }\over {\partial v_{2}}},
 {{\partial\ \ }\over {\partial v_{3}}}, {{\partial\ \ }\over {\partial v_{4}}}\right).
\end{equation*}
Let
$$D_i=\textrm{tr} \left( \left(2\,Y {{\partial\ \ }\over {\partial Y}}\right)^i\right),
\quad i=1,2,3,4$$
and
\begin{equation*}
\Omega_j={{\partial\ \ } \over {\partial V}} \left( 2\,Y {{\partial\ \ }\over {\partial Y}}\right)^j
Y\,  {}^{{}^{{}^{{}^\text{\scriptsize $t$}}}}\!\!\!\left( {{\partial\ \ } \over {\partial V}}\right),
\quad j=0,1,2,3.
\end{equation*}
It is easily seen that
$$ D_1 = \textrm{tr} \left( 2\,Y {{\partial\ \ }\over {\partial Y}}\right) = 2\,\sum_{i=1}^{10}
\, y_i {{\partial\ \ }\over {\partial y_i}}.$$
For brevity we put
$$\partial_{i}= {{\partial\ \ }\over {\partial v_{i} }},\quad i=1,2,3,4.$$
Then we get
\begin{eqnarray*}
\Omega_0&=& y_1\,\partial_1^2 +y_2\,\partial_2^2 + y_3\,\partial_3^2 +y_4\,\partial_4^2 +
2\,y_5\,\partial_1\partial_2 \\
& & + 2\,y_6\,\partial_1\partial_3 + 2\,y_7\,\partial_1\partial_4
+ 2\,y_8\,\partial_2\partial_3 + 2\,y_9\,\partial_2\partial_4 + 2\,y_{10}\,\partial_3\partial_4.
\end{eqnarray*}
\noindent
We observe that $D_1,\,D_2,\,D_3,\,D_4,\,\Omega_0,\,\Omega_1,\Omega_2,\Omega_3$ are invariant differential operators in
$\BD \big( {\mathcal P}_{4,1} \big).$ Then we have the following relations
\begin{equation}
[D_i,D_j]=0 \quad \textrm{for all}\ i,j=1,2,3,4
\end{equation}
and
\begin{equation}
[D_1,\Omega_0]=2\,\Omega_0 .
\end{equation}
Therefore $\BD ({\mathcal P}_{4,1})$ is not commutative.

\vskip 0.87cm


\vskip 0.5cm
\noindent
{\bf 6.2.} {\bf The Case $n=4$ and $m=2$}
\vskip 0.3cm
In this case,
$$GL_{4,2}=GL(4,\BR)\ltimes
\BR^{(2,4)},\ \ K=O(4)\ \ \text{and} \ \ {\mathcal P}_{4,2}=GL_{4,2}/K={\mathcal
P}_4\times \BR^{(2,4)}.$$ We see easily that
$$\fp_{\star}=\,\left\{\,(X,Z)\ \vert\ X=\,^tX\in \BR^{(4,4)},\ \ Z\in
\BR^{(2,4)} \,\right\}.$$

\noindent We put
\begin{eqnarray*}
E_1&=&\begin{pmatrix} 1 & 0 & 0 & 0\\
0 & 0 & 0  & 0\\ 0 & 0 & 0  & 0 \\ 0 & 0 & 0  & 0 \end{pmatrix},\
E_2=\begin{pmatrix} 0 & 0 & 0 & 0\\
0 & 1 & 0  & 0\\ 0 & 0 & 0  & 0 \\ 0 & 0 & 0  & 0 \end{pmatrix},\
E_3=\begin{pmatrix} 0 & 0 & 0 & 0\\
0 & 0 & 0  & 0\\ 0 & 0 & 1  & 0 \\ 0 & 0 & 0  & 0 \end{pmatrix},\\
E_4&=&\begin{pmatrix} 0 & 0 & 0 & 0\\
0 & 0 & 0  & 0\\ 0 & 0 & 0  & 0 \\ 0 & 0 & 0  & 1 \end{pmatrix},
\ E_5=\begin{pmatrix} 0 & {\frac 12} & 0 & 0\\
{\frac 12} & 0 & 0  & 0\\ 0 & 0 & 0  & 0 \\ 0 & 0 & 0  & 0 \end{pmatrix},
\ \ E_6=\begin{pmatrix} 0 & 0 & {\frac 12} & 0\\
0 & 0 & 0  & 0\\ {\frac 12} & 0 & 0  & 0 \\ 0 & 0 & 0  & 0 \end{pmatrix},
\
\\
E_7&=&\begin{pmatrix} 0 & 0 & 0 & {\frac 12} \\
0 & 0 & 0  & 0\\ 0 & 0 & 0  & 0 \\ {\frac 12} & 0 & 0  & 0 \end{pmatrix},\
E_8=\begin{pmatrix} 0 & 0 & 0 & 0 \\
0 & 0 & {\frac 12}  & 0\\ 0 &  {\frac 12}& 0  & 0 \\ 0 & 0 & 0  & 0 \end{pmatrix},\\
E_9&=& \begin{pmatrix} 0 & 0 & 0 & 0 \\
0 & 0 & 0 & {\frac 12}  \\ 0 & 0 & 0  & 0 \\ 0 & {\frac 12} & 0  & 0 \end{pmatrix},\
E_{10}=\begin{pmatrix} 0 & 0 & 0 & 0 \\
0 & 0 & 0 & 0  \\ 0 & 0 & 0  & {\frac 12}  \\ 0 & 0 & {\frac 12} & 0   \end{pmatrix}.
\end{eqnarray*}

Let $O_4$ be the $4\times 4$ zero matrix and let
$$O_{2,4}=\begin{pmatrix} 0 & 0 & 0 & 0 \\ 0 & 0 & 0 & 0 \end{pmatrix}
\in\BR^{(2,4)}$$.
We put
\begin{eqnarray*}
e_i&=& (E_i,O_{2,4}),\quad 1\leq i\leq 10,\\
f_1&=& \left( O_4, \begin{pmatrix} 1 & 0 & 0 & 0 \\ 0 & 0 & 0 & 0 \end{pmatrix}\right),\ f_2=
\left( O_4, \begin{pmatrix} 0 & 1 & 0 & 0 \\ 0 & 0 & 0 & 0 \end{pmatrix}\right),\\
f_3&=& \left( O_4, \begin{pmatrix} 0 & 0 & 1 & 0 \\ 0 & 0 & 0 & 0 \end{pmatrix}\right),
\ f_4= \left( O_4, \begin{pmatrix} 0 & 0 & 0 & 1 \\ 0 & 0 & 0 & 0 \end{pmatrix}\right),\\
f_5&=& \left( O_4, \begin{pmatrix} 0 & 0 & 0 & 0 \\ 1 & 0 & 0 & 0 \end{pmatrix}\right),\ f_6=
\left( O_4, \begin{pmatrix} 0 & 0 & 0 & 0 \\ 0 & 1 & 0 & 0 \end{pmatrix}\right),\\
f_7&=& \left( O_4, \begin{pmatrix} 0 & 0 & 0 & 0 \\ 0 & 0 & 1 & 0 \end{pmatrix}\right),
\ f_8= \left( O_4, \begin{pmatrix} 0 & 0 & 0 & 0 \\ 0 & 0 & 0 & 1 \end{pmatrix}\right).
\end{eqnarray*}
Then
$\{\,e_i,\,f_j\,|\,1\leq i\leq 10,\ 1\leq j\leq 8\,\}$ forms a basis for
$\fp_{\star}.$
Using this basis,
We write for variables $(X,Z)\in \fp_{\star}$ by
$$X=\begin{pmatrix} x_1 & {\frac 12}x_5 & {\frac 12}x_6  & {\frac 12}x_7\\
{\frac 12}x_5 & x_2 & {\frac 12}x_8  & {\frac 12}x_9\\
{\frac 12}x_6 & {\frac 12}x_8 &x_3 & {\frac 12}x_{10}\\
{\frac 12}x_7 & {\frac 12}x_9 & {\frac 12}x_{10} & x_4
\end{pmatrix}
 \ \ \text{and} \ \
Z=\begin{pmatrix} z_{11} & z_{12} & z_{13} & z_{14} \\
z_{21} & z_{22} & z_{23} & z_{24} \end{pmatrix}.$$

We set
\begin{eqnarray*}
\square_{11}&=& {\frac 12}\,\big( 2 A\,x_1 + E\,x_5 + F\,x_6 + G\,x_7 \big),\\
\square_{22}&=&{\frac 12}\,\big( 2 B\, x_2 + E\, x_5 + H\, x_8 + I\, x_9 \big), \\
\square_{33}&=& {\frac 12}\,\big( 2 C\, x_3 + F\, x_6 + H\, x_8 + J\, x_{10} \big), \\
\square_{44}&=&  {\frac 12}\,\big( 2 D\, x_4 + G\, x_7 + I\, x_9 + J\, x_{10} \big), \\
\square_{12}&=& {\frac 12}\,\left\{ 2E(x_1+x_2)+(A+B)x_5+ H x_6 + I x_7 + F x_8 + G x_9 \right\},  \\
\square_{13}&=& {\frac 12}\,\left\{ 2F(x_1+x_3)+ H x_5 + (A+C)x_6+  J x_7 + E x_8 + G x_{10} \right\}, \\
\square_{14}&=& {\frac 12}\,\left\{ 2G(x_1+x_4)+ I x_5 + J x_6 + (A+D)x_7 +  E x_9 + F x_{10} \right\}, \\
\square_{23}&=& {\frac 12}\,\left\{ 2H(x_2+x_3)+ F x_5 + E x_6 + (B+C)x_8 +  J x_9 + I x_{10} \right\}, \\
\square_{24}&=& {\frac 12}\,\left\{ 2I(x_2+x_4)+ G x_5 + E x_7 + J x_8 + (B+D)x_9 +  H x_{10} \right\}, \\
\square_{34}&=& {\frac 12}\,\left\{ 2J(x_3+x_4)+ G x_6 + F x_7 + I x_8 + H x_9 + (C+D) x_{10} \right\}.
\end{eqnarray*}

According to Theorem 3.3, the algebra $\textrm{Pol}(\fp_{\star})^K$ is generated by
the following 16 polynomials

\begin{eqnarray*}
\alpha_1(X,Z)&=& x_1+x_2+x_3+x_4,\\
\alpha_2(X,Z)&=& x_1^2+x_2^2+x_3^2+x_4^2+{\frac 12}\,\big(x_5^2+x_6^2+x_7^2+x_8^2+x_9^2+x_{10}^2\big),\\
\alpha_3(X,Z)&=& x_1^3+x_2^3+x_3^3+x_4^3\\
& & + {\frac 34}\,x_1 \big( x_5^2+x_6^2+x_7^2\big)+ {\frac 34}\,x_2 \big( x_5^2+x_8^2+x_9^2\big)  \\
& & + {\frac 34}\,x_3 \big( x_6^2+x_8^2+x_{10}^2\big)+ {\frac 34}\,x_4 \big( x_7^2+x_9^2+x_{10}^2\big)  \\
& & + {\frac 34}\,\big( x_5\, x_6\, x_8 + x_5 \, x_7\, x_9 + x_6\, x_7\, x_{10} + x_8\, x_9 \, x_{10} \big),\\
\alpha_4(X,Z)&=& A^2+B^2+C^2+D^2+ 2\,\big( E^2+F^2+G^2+H^2+I^2+J^2\big),\\
\beta_{11}^{(0)}(X,Z)&=& z_{11}^2 + z_{12}^2 + z_{13}^2 + z_{14}^2, \\
\beta_{12}^{(0)}(X,Z)&=& z_{11}z_{21} + z_{12} z_{22} + z_{13} z_{23} + z_{14} z_{24}, \\
\beta_{22}^{(0)}(X,Z)&=&  z_{21}^2 + z_{22}^2 + z_{23}^2 + z_{24}^2, \\
\end{eqnarray*}

\begin{eqnarray*}
\beta_{11}^{(1)}(X,Z)&=& x_1 z_{11}^2 + x_2 z_{12}^2 + x_3 z_{13}^2 + x_4 z_{14}^2  + x_5 z_{11} z_{12} \\
& & + x_6 z_{11} z_{13} + x_7 z_{11} z_{14} + x_8 z_{12} z_{13} + x_9 z_{12} z_{14} + x_{10} z_{13} z_{14} ,\\
\beta_{12}^{(1)}(X,Z)&=& x_1 z_{11}z_{21} + x_2 z_{12}z_{22} + x_3 z_{13}z_{23} + x_4 z_{14}z_{24}   \\
& & + {\frac 12}\, x_5\,(z_{11}z_{22}+z_{12}z_{21}) + {\frac 12}\,x_6\,( z_{11}z_{23} + z_{13}z_{21})\\
& & + {\frac 12}\, x_7\,(z_{11}z_{24}+z_{14}z_{21}) + {\frac 12}\,x_8\,( z_{12}z_{23} + z_{13}z_{22})\\
& & + {\frac 12}\, x_9\,(z_{12}z_{24}+z_{14}z_{22}) + {\frac 12}\,x_{10}\,( z_{13}z_{24} + z_{14}z_{23}),\\
\beta_{22}^{(1)}(X,Z)&=& x_1z_{21}^2 + x_2 z_{22}^2 + x_3 z_{23}^2 + x_4 z_{24}^2 + x_5 z_{21} z_{22}  \\
& & + x_6 z_{21} z_{23} + x_7 z_{21} z_{23} + + x_8 z_{22} z_{23} + x_9 z_{22} z_{24} + x_{10} z_{23} z_{24},\\
\beta_{11}^{(2)}(X,Z)&=& A\,z_{11}^2 + B\,z_{12}^2 + C\,z_{13}^2 + D\,z_{14}^2 + 2\,E\,z_{11}z_{12}
+ 2\,F\,z_{11}z_{13} \\
& & +\, 2\,G\,z_{11}z_{14}
+ 2\,H\,z_{12}z_{13} + 2\,I\,z_{12}z_{14}
+ 2\,J\,z_{13}z_{14},\\
\beta_{12}^{(2)}(X,Z)&=& A\,z_{11} z_{21} + B\,z_{12} z_{22} + C\,z_{13} z_{23} + D\,z_{14} z_{24}    \\
& & + \,E\,\big( z_{11}z_{22} + z_{12}z_{21} \big) + \,F\,\big( z_{11}z_{23} + z_{13}z_{21} \big) +
\,G\,\big( z_{11}z_{24} + z_{14}z_{21} \big)    \\
& & + \,H\,\big( z_{12}z_{23} + z_{13}z_{22} \big) + \,I\,\big( z_{12}z_{24} + z_{14}z_{22} \big) +
\,J\,\big( z_{13}z_{24} + z_{14}z_{23} \big),    \\
\beta_{22}^{(2)}(X,Z)&=& A\,z_{21}^2 + B\,z_{22}^2 + C\,z_{23}^2 + D\,z_{24}^2 + 2\,E\,z_{21}z_{22}
+ 2\,F\,z_{21}z_{23} \\
& & +\, 2\,G\,z_{21}z_{24}
+ 2\,H\,z_{22}z_{23} + 2\,I\,z_{22}z_{24}
+ 2\,J\,z_{23}z_{24},  \\
\beta_{11}^{(3)}(X,Z)&=&  \square_{11}\,z_{11}^2 + \square_{22}\, z_{12}^2 +
\square_{33}\,z_{13}^2 + \square_{44}\, z_{14}^2 + \square_{12}\,z_{11}z_{12}   \\
& & +\, \square_{13}\,z_{11}z_{13} + \square_{14}\,z_{11}z_{14} + \square_{23}\,z_{12}z_{13} + \square_{24}\,z_{12}z_{14}
+ \square_{34}\,z_{13}z_{14}, \\
\beta_{12}^{(3)}(X,Z)&=& \square_{11}\,z_{11} z_{21} + \square_{22}\, z_{12} z_{22} +
\square_{33}\,z_{13} z_{23} + \square_{44}\, z_{14} z_{24} + \square_{12}\,z_{11}z_{22}   \\
& & +\, \square_{13}\,z_{11}z_{23} + \square_{14}\,z_{11}z_{24} + \square_{23}\,z_{12}z_{23} + \square_{24}\,z_{12}z_{24}
+ \square_{34}\,z_{13}z_{24},  \\
\beta_{22}^{(3)}(X,Z)&=& \square_{11}\,z_{21}^2 + \square_{22}\, z_{22}^2 +
\square_{33}\,z_{23}^2 + \square_{44}\, z_{24}^2 + \square_{12}\,z_{21}z_{22}   \\
& & +\, \square_{13}\,z_{21}z_{23} + \square_{14}\,z_{21}z_{24} + \square_{23}\,z_{22}z_{23} + \square_{24}\,z_{22}z_{24}
+ \square_{34}\,z_{23}z_{24}.
\end{eqnarray*}
Here $A,B,C,\cdots,J$ are defined in (7.1)-(7.10).

\vskip 0.5cm
We set
\begin{equation*}
\Delta_{ab}:=\det
\begin{pmatrix}
\beta_{11}^{(a)} & \beta_{12}^{(b)} \\
\beta_{12}^{(a)} & \beta_{22}^{(b)}
\end{pmatrix}
\qquad \textrm{for}\ a,b=0,1,2,3.
\end{equation*}

Minoru Itoh claims that
\begin{eqnarray}
&& \big( \alpha_1^3- 3\,\alpha_1\alpha_2\,+\,2\,\alpha_3\big)\,\Delta_{00}\,
 - 3\,\big( \alpha_1^2-\alpha_2\big) \,\big( \Delta_{01}\,+\,\Delta_{10}\big)\\
&&\,+\,6\,\alpha_1\,\big( \Delta_{02}\,+\,\Delta_{11}\,+\,\Delta_{20}\big)
 \,+\,6\,\big(\Delta_{03}\,+\, \Delta_{12}\,+\,\Delta_{21}\,+\,\Delta_{30}\big)=0.\nonumber
\end{eqnarray}
Furthermore he says that there are some more other relations among the above 16 invariants.

\vskip 0.5cm We take a coordinate $(Y,V)$ in ${\mathcal P}_{4,2}$ that is
\begin{equation*}
Y=\begin{pmatrix} y_1 & y_5 & y_6 & y_7 \\ y_5 & y_2 & y_8 & y_9 \\
y_6 & y_8 & y_3 & y_{10} \\ y_7 & y_9 & y_{10} & y_4
\end{pmatrix}
\quad \textrm{and}\quad
V=\begin{pmatrix}
v_{11} & v_{12} & v_{13} & v_{14} \\ v_{21} & v_{22} & v_{23} & v_{24}
\end{pmatrix}.
\end{equation*}

We put
\begin{equation*}
{{\partial\ \ }\over {\partial Y}}=
\begin{pmatrix} {{\partial\ \ }\over {\partial y_1}} & {\frac 12}{{\partial\ \ }\over {\partial y_5}}
& {\frac 12}{{\partial\ \ }\over {\partial y_6}}  & {\frac 12}{{\partial\ \ }\over {\partial y_7}} \\
{\frac 12}{{\partial\ \ }\over {\partial y_5}} & {{\partial\ \ }\over {\partial y_2}}
& {\frac 12}{{\partial\ \ }\over {\partial y_8}} & {\frac 12}{{\partial\ \ }\over {\partial y_9}} \\
{\frac 12}{{\partial\ \ }\over {\partial y_6}} & {\frac 12}{{\partial\ \ }\over {\partial y_8}}
& {{\partial\ \ }\over {\partial y_3}} & {\frac 12}{{\partial\ \ }\over {\partial y_{10}}} \\
{\frac 12}{{\partial\ \ }\over {\partial y_7}} & {\frac 12}{{\partial\ \ }\over {\partial y_9}}
& {\frac 12} { {\partial\ \ }\over {\partial y_{10}} } & { {\partial\ \ }\over {\partial y_4} }
\end{pmatrix}
\quad \textrm{and}\quad
{{\partial\ \ } \over {\partial V}}=
\begin{pmatrix}
{{\partial\quad}\over {\partial v_{11}}} & {{\partial\quad}\over {\partial v_{12}}}
& {{\partial\quad}\over {\partial v_{13}}} & {{\partial\quad}\over {\partial v_{14}}} \\
{{\partial\quad}\over {\partial v_{21}}} & {{\partial\quad}\over {\partial v_{22}}}
& {{\partial\quad}\over {\partial v_{23}}} & {{\partial\quad}\over {\partial v_{24}}}
\end{pmatrix}.
\end{equation*}
Let
\begin{equation*}
D_i=\,\textrm{tr} \left(\left( 2\,Y {{\partial\ \ }\over {\partial Y}} \right)^i\right),\quad i=1,2,3,4
\end{equation*}
and
\begin{equation*}
\Omega_{pq}^{(k)}=\,\left\{\,{{\partial}\over{\partial V}}\left(
2Y{{\partial\ \ }\over{\partial Y}}\right)^kY\,{}^{{}^{{}^{{}^\text{\scriptsize $t$}}}}\!\!\!\left({{\partial\ \ }\over
{\partial V}}\right)\,\right\}_{pq},\quad k=0,1,2,3, \ 1\leq p\leq q\leq 2.
\end{equation*}
We note that $D_1,\,D_2,\,D_3,\,D_4,\,\Omega_{11}^{(0)},\cdots ,\,\Omega_{22}^{(3)}$ are $GL_{4,2}$-invariant.
It is easily seen that
$$ D_1 = \textrm{tr} \left( 2\,Y {{\partial\ \ }\over {\partial Y}}\right) = 2\,\sum_{i=1}^{10}
\, y_i {{\partial\ \ }\over {\partial y_i}}.$$
For brevity we put
$$\partial_{ij}= {{\partial\quad}\over {\partial v_{ij} }},\quad i=1,2,\ 1\leq j\leq 4.$$
Then we get

\begin{eqnarray*}
\Omega_{11}^{(0)} &=& y_1\, \partial_{11}^2 + y_2 \,\partial_{12}^2 + y_3\, \partial_{13}^2 + y_4\, \partial_{14}^2
+ 2\,y_5\, \partial_{11}\partial_{12} +  2\,y_6\, \partial_{11}\partial_{13} \\
& & +    2\,y_7\, \partial_{11}\partial_{14}
+ 2\,y_8\, \partial_{12}\partial_{13} +  2\,y_9\, \partial_{12}\partial_{14} + 2\,y_{10}\, \partial_{13}\partial_{14},\\
\Omega_{12}^{(0)} &=& y_1\,\partial_{11}\partial_{21} + y_2\,\partial_{12}\partial_{22}
+ y_3\,\partial_{13}\partial_{23} + y_4\,\partial_{14}\partial_{24} \\
& & +\, y_5\,\big( \partial_{11}\partial_{22} + \partial_{12}\partial_{21}\big)
+ \, y_6\,\big( \partial_{11}\partial_{23} + \partial_{13}\partial_{21}\big) \\
& & +\, y_7\,\big( \partial_{11}\partial_{24} + \partial_{14}\partial_{21}\big)
+ \, y_8\,\big( \partial_{12}\partial_{23} + \partial_{13}\partial_{22}\big) \\
& & +\, y_9\,\big( \partial_{12}\partial_{24} + \partial_{14}\partial_{22}\big)
+ \, y_{10}\,\big( \partial_{13}\partial_{24} + \partial_{14}\partial_{23}\big), \\
\Omega_{22}^{(0)} &=& y_1\, \partial_{21}^2 + y_2 \,\partial_{22}^2 + y_3\, \partial_{23}^2 + y_4\, \partial_{24}^2
+ 2\,y_5\, \partial_{21}\partial_{22} +  2\,y_6\, \partial_{21}\partial_{23} \\
& & +    2\,y_7\, \partial_{21}\partial_{24}
+ 2\,y_8\, \partial_{22}\partial_{23} +  2\,y_9\, \partial_{22}\partial_{24} + 2\,y_{10}\, \partial_{23}\partial_{24}.
\end{eqnarray*}
Then we have the following relations
\begin{equation}
[D_i,D_j]=0 \quad \textrm{for all}\ i,j=1,2,3,4
\end{equation}
\begin{equation}
[\Omega_{kl}^{(0)},\Omega_{pq}^{(0)}]=0,\quad 1\leq k\leq l\leq 2,\ 1\leq p\leq q\leq 2
\end{equation}
and
\begin{equation}
[D_1,\Omega_{11}^{(0)}]=2\,\Omega_{11}^{(0)},\quad [D_1,\Omega_{12}^{(0)}]=2\,\Omega_{12}^{(0)},\quad
[D_1,\Omega_{22}^{(0)}]=2\,\Omega_{22}^{(0)}.
\end{equation}
Therefore $\BD ({\mathcal P}_{4,2})$ is not commutative.

\end{section}

\vskip 0.87cm


\begin{section}{{\bf Final Remarks}}
\setcounter{equation}{0}
\vskip 0.3cm
In this section we present some open problems and discuss a notion of automrphic forms on ${\mathcal P}_{n,m}$
using $GL_{n,m}$-invariant differential operators on the Minkowski-Euclid space ${\mathcal P}_{n,m}.$
\vskip 0.3cm
We recall the invariant polynomials $\alpha_j\, (1\leq j\leq n)$ in (3.11) and
$\beta_{pq}^{(k)}\ (0\leq k\leq n-1,\ 1\leq p\leq q\leq m)$ in (3.12). We also recall the invariant differential
operators $D_j\, (1\leq j\leq n)$ in (3.19) and
$\Omega_{pq}^{(k)}\ (0\leq k\leq n-1,\ 1\leq p\leq q\leq m)$ in (3.20).

\vskip 0.5cm\noindent
\begin{theorem}
The following relations hold:
\begin{equation}
[D_i,D_j]=0\quad \textrm{for all}\ 1\leq i,j\leq n,
\end{equation}
\begin{equation}
[\Omega_{kl}^{(0)},\Omega_{pq}^{(0)}]=0,\quad 1\leq k\leq l\leq m,\ 1\leq p\leq q\leq m
\end{equation}
and
\begin{equation}
[D_1,\Omega_{pq}^{(0)}]=2\,\Omega_{pq}^{(0)}\quad \textrm{for all}\ \ 1\leq p\leq q \leq m.
\end{equation}
\end{theorem}
\vskip 0.3cm\noindent
{\it Proof.} The relation (8.1) follows from the work of Atle Selberg (cf.\,\cite{M2, S1,T}).
We take a coordinate $(Y,V)$ in ${\mathcal P}_{n,m}$ with $Y=(y_{ij})$ and $V=(v_{kl})$.
We put
$$\Yd=\left( { {1+\delta_{ij}}\over 2 }
{{\partial}\over {\partial y_{ij} }}\right)\quad \textrm{and}\quad \Vd=\left({
{\partial}\over {\partial v_{kl}} } \right),$$ where $1\leq
i,j,l\leq n$ and $1\leq k\leq m.$ Then we get

\begin{eqnarray*}
D_1 &=& 2\,\sum_{1\le i\leq j\leq n} y_{ij}\,
{{\partial}\over {\partial y_{ij} }},\\
\Omega_{pq}^{(0)}&=& \sum_{a=1}^n y_{aa}\,{{\partial^2\qquad}\over {\partial v_{pa}\partial v_{qa}} }\,+\,
\sum_{1\leq a < b\leq n} y_{ab} \left( {{\partial^2\qquad}\over {\partial v_{pa}\partial v_{qb}} }\,+\,
{{\partial^2\qquad}\over {\partial v_{pb}\partial v_{qa}} } \right).
\end{eqnarray*}
By a direct calculation, we obtain the desired relations (8.2) and (8.3).
\hfill $\square$

\vskip 0.5cm
\noindent
{\bf Conjecture 2.}
\begin{eqnarray}
\Theta_{n,m}(\alpha_j)&=& D_j\quad \textrm{for all}\ 1\leq j\leq n,\\
\Theta_{n,m}\big(\beta_{pq}^{(k)}\big)&=& \Omega_{pq}^{(k)}\quad \textrm{for all}\ 0\leq k\leq n-1,\ 1\leq p\leq q\leq m.
\end{eqnarray}

\vskip 0.3cm We refer to Conjecture 1 in Section 2.

\vskip 0.5cm
\noindent
{\bf Conjecture 3.} The invariants $D_j\, (1\leq j\leq n)$and
$\Omega_{pq}^{(k)}\ (0\leq k\leq n-1,\ 1\leq p\leq q\leq m)$ generate the noncommutative algebra
$\BD ({\mathcal P}_{n,m})$.

\vskip 0.5cm
\noindent
{\bf Conjecture 4.} The above relations (8.1), (8.2) and (8.3)
generate all the relations among the set
$$\left\{ D_j,\,\Omega_{pq}^{(k)}\,|\
1\leq j\leq n,\ 0\leq k\leq n-1,\ 1\leq p\leq q\leq m \right\}.$$

\vskip 0.5cm \noindent
{\bf Problem 8.} Find a natural way to construct generators of $\BD ({\mathcal P}_{n,m})$.

\vskip 0.5cm
Using $GL_{n,m}$-invariant differential operators on the Minkowski-Euclid space ${\mathcal P}_{n,m},$
we introduce a notion of automorphic forms on ${\mathcal P}_{n,m}.$
\vskip 0.3cm\noindent
Let
\begin{equation*}
\Gamma_{n,m}:=GL(n,\BZ)\ltimes \BZ^{(m,n)}
\end{equation*}
be the arithmetic subgroup of $GL_{n,m}.$
Let ${\mathcal Z}_{n,m}$ be the center of $\BD ({\mathcal P}_{n,m})$.

\vskip 0.5cm\noindent
\begin{definition} A smooth function $f:{\mathcal P}_{n,m}\lrt \BC$ is said to be an automorphic form
for $\Gamma_{n,m}$ if it satisfies the following conditions:
\vskip 0.2cm\noindent
(A1) $f$ is $\Gamma_{n,m}$-invariant.
\vskip 0.2cm\noindent
(A2) $f$ is an eigenfunction of any differential operator in the center ${\mathcal Z}_{n,m}$
of $\BD ({\mathcal P}_{n,m})$.
\vskip 0.2cm\noindent
(A3) $f$ has a growth condition.
\end{definition}

\vskip 0.5cm \noindent
{\bf Problem 9.} Find the center ${\mathcal Z}_{n,m}$
of $\BD ({\mathcal P}_{n,m})$ explicitly.
\vskip 0.5cm\noindent
{\bf Problem 10.} Find the center of the universal enveloping algebra of the complexification of the Lie algebra
of the group $GL_{n,m}$ explicitly.

\vskip 0.5cm
We may define another notion of automorphic forms as follows.

\begin{definition} Let $\BD_\spadesuit$ be a commutative
subalgebra of $\BD ({\mathcal P}_{n,m})$ containing the Laplacian $\Delta_{n,m;A,B}.$
A smooth function $f:{\mathcal P}_{n,m}\lrt \BC$ is said to be an automorphic form
for $\Gamma_{n,m}$ with respect to $\BD_\spadesuit$ if it satisfies the following conditions:
\vskip 0.2cm\noindent
(A1) $f$ is $\Gamma_{n,m}$-invariant.
\vskip 0.2cm\noindent
(A2) $f$ is an eigenfunction of any differential operator in $\BD_\spadesuit$.
\vskip 0.2cm\noindent
(A3) $f$ has a growth condition.
\end{definition}

\end{section}



\vskip 1cm
\begin{thebibliography}{99}






\bibitem{HC1} Harish-Chandra, \textit{Representations of semisimple Lie groups,
I.,} Trans. Amer. Math. Soc. {\bf 75} (1953), 185-243.

\bibitem{HC2} Harish-Chandra, \textit{The characters of semisimple Lie groups,
} Trans. Amer. Math. Soc. {\bf 83} (1956), 98-163.


\bibitem{He1} S. Helgason, \textit{Differential operators
on homogeneous spaces,} Acta Math. {\bf 102} (1959), 239-299.

\bibitem{He2} S. Helgason, \textit{Groups and geometric
analysis,} Academic Press, New York (1984).
\bibitem{Ho} R. Howe, \textit{Perspectives on invariant theory: Schur duality,
multiplicity-free actions and beyond,} The Schur lectures
(1992)\,(Tel Aviv), Israel Math. Conf. Proceedings, {\bf vol.\ 8}
(1995), 
1-182.

\bibitem{M1} H. Maass, \textit{Die Bestimmung der Dirichletreihnen mit Gr{\"o}ssencharakteren
zu den Modulformen $n$-ten Grades,} J. Indian Math. Soc., {\bf
vol. 19} (1955), 1-23.


\bibitem{M2} H. Maass, \textit{Siegel modular forms and Dirichlet
series,} Lecture Notes in Math., {\bf vol. 216}, Springer-Verlag,
Berlin-Heidelberg-New York (1971).

\bibitem{Mi} H. Minkowski, \textit{Gesammelte Abhandlungen,}
Chelsea, New York (1967).

\bibitem{S1} A. Selberg, \textit{Harmonic analysis and
discontinuous groups in weakly symmetric Riemannian spaces with
applications to Dirichlet series,} J. Indian Math. Soc. B.  {\bf
20} (1956), 47-87.




\bibitem{T} A. Terras, \textit{Harmonic Analysis on
Symmetric Spaces and Applications II,} Springer-Verlag (1988).

\bibitem{W} H. Weyl, \textit{The classical groups: Their
invariants and representations,} Princeton Univ. Press, Princeton,
New Jersey, second edition (1946).



\bibitem{Y1} J.-H. Yang, \textit{Singular
Jacobi forms,} Trans. of American Math. Soc. {\bf 347}, No. 6
(1995), 2041-2049.

\bibitem{Y2} J.-H. Yang, \textit{Polarized Real Tori,}
arXiv:0912.5084v1 [math.AG] (2009).




\end{thebibliography}
\end{document}